\documentclass[12pt, abstracton]{scrartcl}   

 \usepackage[T1]{fontenc}              
 \usepackage[latin1]{inputenc}        
      \input{dehyphtex}              
 \usepackage[width=16cm, height=22cm]{geometry}   
 \usepackage{amsmath,amssymb}         
 \usepackage[amsmath,thmmarks,hyperref]{ntheorem} 
 \usepackage{icomma}             
 \usepackage{enumerate}             
 \usepackage{mathrsfs}         
 \usepackage{slashed}
 \usepackage{color}

 \numberwithin{equation}{section}
 
 \usepackage[numbers,square]{natbib}
\theoremstyle{nonumberplain}  
\theoremheaderfont{\itshape}  
\theorembodyfont{\normalfont}  
\theoremseparator{.}  
\theoremsymbol{$\Box$}
\newtheorem{proof}{Proof} 
  
\theoremstyle{plain}  
\theoremheaderfont{\normalfont\bfseries}  
\theorembodyfont{\itshape}  
\theoremsymbol{~}
\theoremseparator{.}  
\newtheorem{proposition}{Proposition}[section]  
\newtheorem{corollary}[proposition]{Corollary}  
\newtheorem{lemma}[proposition]{Lemma}  
\newtheorem{theorem}[proposition]{Theorem}   

\theorembodyfont{\normalfont}  
 
\newtheorem{remark}[proposition]{Remark}
\newtheorem{example}[proposition]{Example}  
  
\newtheorem{definition}[proposition]{Definition} 

\theoremstyle{nonumberplain}
\theorembodyfont{\normalfont}

\newcommand*{\grad}{\operatorname{grad}}

\newcommand{\R}{\mathbb{R}}
\newcommand{\Eh}{\mathcal{E}}

\newcommand{\N}{\mathbb{N}}
\newcommand{\Z}{\mathbb{Z}}

\newcommand{\dd}{\mathrm{d}}

\newcommand{\ind}{\operatorname{ind}}
\newcommand{\Cl}{\mathcal{C}\!\ell}
\newcommand{\tr}{\mathrm{tr}}
\newcommand{\End}{\mathrm{End}}

\renewcommand*{\div}{\operatorname{div}}
\newcommand{\cc}{\mathbf{c}}
\newcommand{\bb}{\mathbf{b}}
\DeclareMathOperator{\str}{\mathrm{str}}
\newcommand{\vol}{\mathrm{vol}}

\newcommand{\TMd}{T^*\!M}

\newcommand{\stimes}{\,\widehat{\otimes}\,}
\renewcommand{\hat}{\widehat}

\title{A Semiclassical Heat Kernel Proof of the Poincaré-Hopf Theorem}
\author{ Matthias Ludewig}

\begin{document}

\maketitle

\begin{center}
  Universität Potsdam / Institut für Mathematik \\ 
  Am Neuen Palais 10 / 14469 Potsdam, Germany \\ \medskip
 matthias.ludewig@uni-potsdam.de \\ \medskip
 Tel.\: +49-331-977-1248
\end{center}

\begin{abstract}
We consider the semiclassical asymptotic expansion of the heat kernel coming from Witten's perturbation of the de Rham complex by a given function. For the index, one obtains a time-dependent integral formula which is evaluated by the method of stationary phase to derive the Poincaré-Hopf theorem. We show how this method is related to approaches using the Thom Form of Mathai and Quillen. Afterwards, we use a more general version of the stationary phase approximation in the case that the perturbing function has critical submanifolds to derive a degenerate version of the Poincaré-Hopf theorem.
\end{abstract}


\section{Introduction}

\dictum[David Hilbert]{Wenn man 2 Wege hat, so muss man nicht bloss diese Wege gehen oder neue suchen, sondern dann das ganze zwischen den beiden Wegen liegende Gebiet erforschen. 

\medskip

\itshape{Given two routes, it is not right to take either of these two or to look for a third; it is necessary to investigate the area lying between the two routes.} \footnotemark}
\footnotetext{Cited in \cite{hughes}}
\bigskip

The connection between semiclassical Analysis and Morse theory was discovered by Witten \cite{witten} about thirty years ago and got a lot of attention ever since. The main idea is inspired by heuristics coming from physics and goes as follows: The Laplace part of a Schrödinger operator $\hbar^2 \Delta + V$ is responsible for diffusion, while the potential term $V$ causes concentration at its cavities. Therefore, when simultaneously taking limits $\hbar \downarrow 0$ and $t \uparrow \infty$ in some appropriate way (where $t$ is the time parameter in the heat equation corresponding to the operator $\hbar^2\Delta + V$), one expects the particles to concentrate near the minima of the potential. With Morse theory on the other hand, one can obtain topological information of the underlying space by investigating the critical points of a given function.

Witten intertwines the Euler operator acting on differential forms with a vector field $X$, this perturbation depending on a small parameter $\hbar$. That way he obtains a Schrödinger type operator whose semiclassical eigenfunction approximations he uses to construct the Morse complex.

The Morse inequalities directly imply the Poincaré-Hopf theorem, which states that the Euler characteristic of $M$ is determined by the critical points of $X$, more precisely
\begin{equation*}
  \chi(M) = \sum_{\{X(p) = 0\}} (-1)^{\nu(p)},
\end{equation*}
where the index $\nu(p)$ of a critical point $p$ is equal to the number of negative eigenvalues of the linearization $\nabla X|_p$. 

Another way to get this theorem uses the Thom form $U$ of Mathai and Quillen: By the transgression formula for the Thom form, the pullbacks of $U$ along any two vector fields are cohomologous; on the other hand, the Euler class of $M$ is the pullback of $U$ along the zero vector field. The pullback of $U$ along the vector field $X_t := t^{1/2} X$ gives a differential form on $M$ and the Poincaré-Hopf theorem follows then by evaluating $X_t^*U$ with the method of stationary phase \citep[Thm.\ 1.56]{bgv}.

These proofs seem conceptually very different at first: Witten uses the low-lying eigenfunctions of his operator to construct a complex chain homotopic to the de Rham complex whence the theorem follows by an argument of homological algebra, while the other proof derives an integral formula that interpolates between the Gauss-Bonnet-Chern theorem and Poincaré-Hopf. In this article, we show that in fact one can use the semiclassical asymptotics \cite{baerpfaeffle} of the heat kernel of Witten's operator to recover the interpolation formula appearing in the Thom form proof.

More precisely, from the semiclassical heat kernel asymptotics, we derive the integral formula 
\begin{equation*}
  \chi(M) =  \int_M \alpha(t)e^{-t|X|^2} ~~~~ \text{for all} ~ t>0,
\end{equation*}
where $\alpha(t)$ is some function depending polynomially on $t$. We then use ideas from \mbox{Getzler's} proof of the local index theorem \citep{getzlerlocal} to explicitly calculate $\alpha(t)$ in terms of the curvature of $M$ and Taylor coefficients of $X$. It turns out that the integrand above is nothing but the Thom form, pulled back via $X_t$.

We show that in the limit $t \downarrow 0$, the integral formula above yields the Gauss-Bonnet-Chern theorem and in the limit $t \uparrow \infty$, the integral can be evaluated with the method of stationary phase to obtain the Poincaré-Hopf theorem. In this sense, the $t$-dependent integral formula above interpolates between two classical theorems.

More generally, by allowing the critical set to be a disjoint union of submanifolds of $M$, we obtain a degenerate version of the Poincaré-Hopf theorem (Thm.\ \ref{ThmMorseBott}). Let us remark that the degenerate Morse inequalities (Morse-Bott inequalities) that imply the degenerate Poincaré-Hopf theorem have been proved with heat kernel methods by Bismut \cite{bismut}, but the integrands are not explicitly calculated in terms of curvature. For other treatments of the Morse-Bott inequalities, see \cite{bottnondegenerate} or \cite{austinbraam}, for example.

This article is organized as follows: At first, we briefly review some basic notions regarding the Clifford algebra and the exterior algebra of a Euclidean vector space $V$. In Section \ref{SectionHeatkernel}, we state the needed results about the semiclassical expansion of the heat kernel of Witten's operator. Afterwards, we introduce Getzler symbols and use them to calculate the relevant integrands. In Section \ref{SectionMcKeanSingr}, we put these results together to prove the stated integral formula for the index and the consequences of it. The last section is dedicated to the degenerate case, in which we derive the Morse-Bott Theorem \ref{ThmMorseBott}. This needs a longer calculation because of the more complicated nature of the corresponding stationary phase expansion.

\medskip

{\bf Acknowledgments.} It is a pleasure to thank Christian Bär and Florian Hanisch for helpful discussion as well as Potsdam Graduate School and the Fulbright Commission for financial support.

\section{The Clifford Symbol on the Exterior Algebra} \label{SectionCliffordSymbol}

We begin with a short recap of the filtrations, gradings and symbol maps associated to the exterior algebra. 

In this section, let $V$ be a Euclidean vector space of dimension $n$. We consider two different gradings on the exterior algebra $\Lambda V$: the $\mathbb{Z}$-grading induced by the degree of forms and the $\mathbb{Z}_2$-grading, where the even (odd) part is the space of even-degree (odd-degree) forms. The latter gives $\Lambda V$ the structure of a superalgebra; the associated grading operator $\Xi$ is by definition the operator that is the identity on even-degree forms and minus the identity on odd-degree forms.

For $v \in V$, we define the Clifford multiplications $\Lambda V \longrightarrow \Lambda V$ by
\begin{equation} \label{CliffordRelatations}
  \cc(v) = \varepsilon(v) - \iota(v), ~~~~~~ \bb(v) = \varepsilon(v) + \iota(v).
\end{equation}
Here, $\varepsilon$ and $\iota$ denote exterior and interior multiplication, where to define the latter, we use the Euclidean structure on $V$. For $v, w \in V$, we have the Clifford relations
\begin{equation*}
  [\cc(v), \cc(w)]_s  = - 2\langle v, w \rangle, ~~~~ [\bb(v), \bb(w)]_s = 2\langle v, w \rangle, ~~~~ [\cc(v), \bb(w)]_s = 0,
\end{equation*}
where $[\,\cdot\, , \,\cdot\,]_s$ denotes the super commutator. Under the isomorphism (of vector spaces) $\Lambda V \cong \Cl(V)$ given by
\begin{equation*}
  v_1 \wedge \dots \wedge v_k \longmapsto v_1 \cdots v_n,
\end{equation*}
$\cc(v)$ acts on a homogeneous element $a \in \Cl(V)$ by Clifford multiplication with $v$ from the left, while $\bb(v)$ acts by Clifford multiplication with $v$ from the right, followed by multiplication with $(-1)^{|a|}$, where $|a|$ is the Clifford order of $a$. The endomorphism space $\End(\Lambda V)$ is generated as an algebra by the elements $\cc(v), \bb(v)$ with $v \in V$. 

\begin{lemma} \label{RepresentationOfGradingOperator}
The grading operator $\Xi$ can be represented by
\begin{equation*}
  \Xi = (-1)^{\frac{n(n+1)}{2}} \, \cc^1 \cdots \cc^n \bb^1 \cdots \bb^n
\end{equation*}
where we wrote $\cc^j = \cc(e^j)$, $\bb^j = \bb(e^j)$ for an orthonormal basis $e^1, \dots, e^n$ of $V$. $\Xi$ does not depend on the choice of this orthonormal basis.
\end{lemma}

\begin{proof}
 Let $I=(i_1, \dots, i_l)$, $i_1 < \dots < i_l$ be a multi-index. Writing $e^I = e^{i_1} \cdots e^{i_l}$, we have
\begin{align*}
  \bb^1 \cdots \bb^n e^I 
  &= (-1)^{l} \bb^1 \cdots \bb^{n-1} e^I \cdot e^n 
  = (-1)^{l+ l+1} \bb^1 \cdots \bb^{n-2} e^I \cdot e^n \cdot e^{n-1}  \\
  &= \dots = (-1)^{nl + 1 + 2 + \dots + {n-1}} e^I \cdot e^n \cdots e^1
  = (-1)^{nl + \frac{n(n-1)}{2}} e^I \cdot e^n \cdots e^1.
\end{align*}
Now
\begin{equation*}
  e^j \cdot e^I \cdot e^j = 
  \begin{cases} 
  (-1)^{l+1}e^I & \text{if} ~~ j \notin I \\
  (-1)^{l}e^I & \text{if} ~~ j \in I
  \end{cases},
\end{equation*}
so that
\begin{align*}
  \cc^1\cdots\cc^n\bb^1\cdots\bb^n e^I 
  = (-1)^{nl + \frac{n(n-1)}{2}} e^1 \cdots e^n \cdot e^I \cdot e^n \cdots e^1 
  = (-1)^{\frac{n(n-1)}{2} + \sum_{j \notin I} 1} \,e^I
\end{align*}
Now  we have $\sum_{j \notin I} 1 = n -l \equiv -n + l \mod 2$ so that
\begin{equation*}
  (-1)^{\frac{n(n+1)}{2}} \cc^1\cdots\cc^n\bb^1\cdots\bb^n e^I = (-1)^l e^I,
\end{equation*}
which was the claim.
\end{proof}

The algebra $\End(\Lambda V)$ has a filtration and a bi-filtration, which we will both use: An element $A \in \End(\Lambda V)$ has (Clifford-) bi-order $(k, l)$ or lower if it can be written as 
\begin{equation} \label{RepresentationOfA}
  A = \sum_{|I| \leq k}\sum_{|J| \leq l} A_{IJ} \,\cc^I \bb^J,
\end{equation}
where we wrote $\cc^I := \cc^{i_1} \cdots \cc^{i_m}$ for a multiindex $I = (i_1 < \dots < i_m)$ and similarly for $\bb^J$. Here $|I|:=m$ is the length of the multi-index. We say that $A$ has (Clifford-) order $k$ or lower if it has bi-order $(k, n)$ or lower.

For such an element $A$ of order $k$, we define its $k$-th Clifford symbol by
\begin{equation*}
  \slashed{\sigma}_k(A) := \slashed{\sigma}_{k, \bullet}(A) := \sum_{|I| = k} \sum_{|J| \leq n} A_{IJ} \, e^I \,\stimes\, \bb^J \in \Lambda V \,\stimes\, \Cl(-V),
\end{equation*}
where $\stimes$ denotes the super tensor product and $\Cl(-V)$ denotes subalgebra\footnote{By \eqref{CliffordRelatations}, this is a Clifford algebra associated to the negative of the scalar product on $V$, hence the notation.} of endomorphisms generated by the $\bb^j$, which super-commutes with the Clifford action induced by $\cc$. The Clifford bi-symbol of an element $A$ of bi-order $(k, l)$ is defined as
\begin{equation*}
  \slashed{\sigma}_{k,l}(A) := \sum_{|I|=k}\sum_{|J|=l} A_{IJ} \, e^I \,\stimes\, e^J \in \Lambda V \,\stimes\, \Lambda V.
\end{equation*}
It is straightforward to check that all these definitions do not depend on the choice of orthonormal basis. Furthermore, both symbols just defined are compatible with the multiplication in the following sense: If $A, B \in \End(V)$ have Clifford order $k$ and $u$ (or bi-order $(k, l)$, $(u, v)$), then $AB$ is of Clifford order $k+u$ (respectively bi-order $(k+u, l+v)$) and we have
\begin{equation} \label{HomomorphismLike}
  \slashed{\sigma}_{k+u}(AB) = \slashed{\sigma}_{k}(A) \slashed{\sigma}_{u}(B) ~~~~~~~~\text{and}~~~~~~
  \slashed{\sigma}_{k+u, l+v}(AB) = \slashed{\sigma}_{k, l}(A) \slashed{\sigma}_{u, v}(B)
\end{equation}
This identity would be wrong if we hadn't used the super tensor product in the target spaces.

\begin{definition}
The supertrace on $\End(\Lambda V)$ is defined by $\str(A) := \tr(\Xi A)$, where $\Xi$ is the grading operator.
\end{definition}

\begin{proposition} \label{SupertraceOfGradingOperator}
If $A$ has bi-order $(k, l)$ with $k <n$ or $l <n$, then $\str\, (A) = 0$. On the other hand, $\str (\Xi) = 2^n$.
\end{proposition}

\begin{proof}
Let $A = \cc^{j_1} \cdots \cc^{j_k} \bb^{i_1} \cdots \bb^{i_l} \in \End(\Lambda V)$ with $k<n$ and let $u \notin \{j_1, \dots j_k\}$. Then straightforward calculation shows
\begin{equation*}
  [ \cc^u \cc^{j_1} \cdots \cc^{j_k} \bb^{i_1} \cdots \bb^{i_l}, \cc^u]_s = (-1)^{k+l} 2 \,\cc^{j_1} \cdots \cc^{j_k} \bb^{i_1} \cdots \bb^{i_l},
\end{equation*}
so that $A$ is a super-commutator. However, $\str$ vanishes on super-commutators \cite[Prop.\ 1.31]{bgv}. In the case that $k = n$ and $l<n$, replace $\cc^u$ by $\bb^u$ for some $u \notin \{i_1, \dots, i_l\}$.

Finally, by definition, $\str(\Xi) = \tr(\Xi^2) = \tr(\mathrm{id}) = 2^n$.
\end{proof}

\begin{corollary} \label{SupertraceFormula}
  We have
  \begin{equation*}
    \str(A) = (-1)^{\frac{n(n+1)}{2}}\,2^n \,\bigl\langle \slashed{\sigma}_{n,n}(A), \vol \, \widehat{\otimes} \, \vol \bigr\rangle
  \end{equation*}
  where $\vol = e^1 \wedge \dots \wedge e^n$ for an (oriented) orthonormal basis $e^1, \dots, e^n$.
\end{corollary}

\begin{remark} \label{RemarkOrientation}
 Notice that this does not depend on an orientation of $M$. Changing the orientation turns $\vol$ into $-\vol$, but the two signs cancel each other.
\end{remark}

\begin{proof}
For $A$ of the form \eqref{RepresentationOfA} (with $k=l=n$), we have
\begin{equation*}
  \str(A) = \str\Bigl(\sum\nolimits_{I, J} A_{IJ} \cc^I \bb^J\Bigr) = \str(A_{[n][n]} \cc^1\cdots \cc^n \bb^1\cdots\bb^n)
  = A_{[n][n]} (-1)^{\frac{n(n+1)}{2}} 2^n,
\end{equation*}
by Lemma \ref{RepresentationOfGradingOperator} and Prop.\ \ref{SupertraceOfGradingOperator}, where we wrote $[n]$ for the multi-index $(1, 2, \dots, n)$.
On the other hand,
\begin{equation*}
  \slashed{\sigma}_{n,n}(A) = A_{[n][n]} \vol \stimes \vol,
\end{equation*}
which gives the result.
\end{proof}


\section{The Semiclassical Heat Kernel Expansion} \label{SectionHeatkernel}

\begin{definition} \label{DefinitionSchroedinger}
Let $(M, g)$ be a Riemannian manifold of dimension $n$ and let $\Eh$ be a real or complex vector bundle  over $M$ equipped with a scalar product (or positive definite Hermitian form respectively). We say that an operator of the form
\begin{equation}
  H_\hbar = \hbar^2 L + \hbar W + V,
\end{equation}
is of Schrödinger type, if $L$ is a formally self-adjoint operator of Laplace type, i.e.\ locally it has the form
\begin{equation*}
  L = - \mathrm{id}_\Eh \sum_{ij =1}^n g^{ij} \frac{\partial^2}{\partial x^i \partial x^j} + \text{lower order terms},
\end{equation*}
 and $V$ and $W$ are symmetric endomorphism fields (we adopt the sign convention such that the eigenvalues of Laplace type operators tend to $+ \infty$).
\end{definition}

\begin{example}[The Witten Operator] \label{ExampleWittenOperator}
 Let $\Eh = \Lambda \TMd$. For a one-form $\xi \in \Omega^1(M)$, define
\begin{equation*}
  D_\hbar = \hbar D + \bb(\xi), ~~~~ \text{where} ~~~~ D = \dd + \dd^*.
\end{equation*}
Then $D_\hbar^2$ is of Schrödinger type, as straightforward calculation shows that in local coordinates, we have
\begin{equation}
  W = [D, \bb(\xi)] = \sum_{i=1}^n \cc(\dd x^i) \, \bb(\nabla_i \xi), ~~~~~ \text{and} ~~~~~ V = |\xi|^2.
\end{equation}
Witten  {\normalfont\citep{witten}} uses this operator, but with a different normalization ($\hbar \mapsto t^{-1}$). Also see {\normalfont \citep[p.\ 125]{johnroe}} and \citep[Chapter 5]{zhang}.
\end{example}

Suppose that $M$ is compact. Then it is well-known that for each $\hbar>0$ fixed, a Schrödinger type operator $H_\hbar$ is an unbounded operator in $L^2(M, \Eh)$ which is self-adjoint on the Sobolev space $H^2(M, \Eh)$ and has eigenvalues tending to $+\infty$ (see e.g.\ \citep[1.6]{gilkey}). For $t>0$, the operator $e^{-tH_\hbar}$ (defined by functional calculus) has a smooth integral kernel
\begin{equation}
  k_\hbar  ~\in~  \Gamma^\infty\bigl(M \times M \times (0, \infty), \Eh^* \boxtimes \Eh\bigr).
\end{equation}
that depends smoothly on $\hbar$ \citep[Thm.\ 2.48]{bgv}. Here $\Eh^* \boxtimes \Eh$ denotes the bundle $\mathrm{pr}_1^* \Eh \otimes \mathrm{pr}_2^* \Eh$ over $M\times M$, where $\mathrm{pr}_{1,2}: M \times M \longrightarrow M$ are projections on the first and second factor respectively. Its fiber over $(p, q) \in M \times M$ is given by $\Eh^*_p \otimes \Eh_q \cong \mathrm{Hom}(\Eh_p, \Eh_q)$. 

On any complete Riemannian manifold $M$, we call
\begin{equation*}
  \mathrm{e}_\hbar(p, q, t) = (4\pi t\hbar^2)^{-n/2} \exp\left( -\frac{1}{4t\hbar^2} d(p, q)^2 \right)
\end{equation*}
the Euclidean heat kernel, because $k_\hbar = \mathrm{e}_\hbar$ if $L$ is the usual Laplace operator on functions in Euclidean space and $W = V = 0$. Here, $d(p, q)$ is the Riemannian distance between $p$ and $q$, so that $\mathrm{e}_\hbar$ is smooth on the set
\begin{equation*}
M \bowtie M = \bigl\{ (p, q) \in M \times M \mid p~~ \text{is not a cut-point of} ~~q \bigr\},
\end{equation*}
which is a dense open neighborhood of the diagonal. (Two points $(p, q)$ are cut-points of each other if either there are several geodesics of minimal length joining $p$ and $q$ or there is a Jacobi field along the unique shortest geodesic connecting $p$ and $q$ that vanishes at both $p$ and $q$ \cite{baerpfaeffle}.)

The $\hbar$-asymptotics of the heat kernel $k_\hbar$ were developed in \citep{baerpfaeffle}. We briefly recall the relevant ideas: Because $k_\hbar$ is the solution kernel to the heat equation, we have $\left(\partial_t + H_\hbar \right) k_\hbar(\,\cdot\,, q) = 0$ for each $q \in M$ and $\hbar>0$. If on $M \bowtie M$ we make the ansatz
\begin{equation}
  k_\hbar \sim \mathrm{e}_\hbar \sum\nolimits_{j=0}^\infty \hbar^j \Phi_j
\end{equation}
for an asymptotic expansion of $k_\hbar$ in $\hbar$, we can formally (i.e.\ termwise) apply the heat operator $(\partial_t + H_\hbar)$  to this expression and straightforward calculation shows that the result is again of the form $\mathrm{e}_\hbar$ times some power series in $\hbar$. One obtains that in order to have this power series vanish, the $\Phi_j$ have to fulfill the recursive transport equations
\begin{equation} \label{RecursiveTransportEquations}
  \left(t\frac{\partial}{\partial t} + \nabla_{\mathcal{V}} + G + tV\right) \Phi_j(\,\cdot\,, q) = - t W \Phi_{j-1}(\,\cdot\,, q) - t L \Phi_{j-2}(\,\cdot\,, q),
\end{equation}
for each $q \in M$, where
\begin{equation} \label{VandG}
  \mathcal{V} = \frac{1}{2}\grad[d(\,\cdot\,, q)^2], ~~~~~~~~~ G = - \left(\frac{n}{2} + \frac{1}{4} \Delta[d(\,\cdot\,, q)^2] \right).
\end{equation}
Therefore one makes the following definition.

\begin{definition}
Let $\sum_{j=0}^\infty \hbar^j \Phi_j$ be a formal power series with coefficients $\Phi_j$ in the space $\Gamma^{\infty}\bigl(M \bowtie M \times [0, \infty), \Eh^* \boxtimes \Eh \bigr)$.
Then the formal expression 
\begin{equation*}
  \hat{k}_\hbar := \mathrm{e}_\hbar \sum\nolimits_{j=0}^\infty \hbar^j \Phi_j
\end{equation*} 
is called semiclassical heat kernel expansion if for all $q \in M$, the coefficients $\Phi_j$ fulfill the recursive transport equations \eqref{RecursiveTransportEquations} with the initial condition $\Phi_0(q, q, 0) = \mathrm{id}_{\Eh_q}$.
\end{definition}

Regarding this, we have the following theorem.

\begin{theorem}[Bär, Pfäffle] \label{ThmBP0} {\normalfont\citep[Lemma 3.1]{baerpfaeffle}}
 Let $M$ be a complete Riemannian manifold, $\Eh$ a vector bundle with scalar product over $M$ and $H_\hbar$ be of Schrödinger type. Then there exists a unique semiclassical  heat kernel expansion. 
 
Furthermore, for each $q \in M$, the first transport equation \eqref{RecursiveTransportEquations} has a unique solution for each prescribed value for $\Phi_0(q, q, 0)$, while for $j>0$, the $j$-th transport equation has a unique solution for each right hand side.
\end{theorem}

For the solution theory of equations of the form \eqref{RecursiveTransportEquations}, also see \cite{ludewigsource}. 

The relation between the formal heat kernel and the true heat kernel is the following.

\begin{theorem}[Bär, Pfäffle] \label{ThmBP} {\normalfont\citep[Thm.\ 3.3]{baerpfaeffle}}
Let $M$ be compact. Then for all $T>0$, $m, k \in \N_0$ and $N > 2n + 2m + 2k$, there exist constants $C>0$, $\hbar_0>0$ such that
\begin{equation*}
  \sup_{t \in (0, T]} \left\| \frac{\partial^k}{\partial t^k} \left( k_\hbar - \chi\mathrm{e}_\hbar\sum\nolimits_{j=0}^N \hbar^j \Phi_j \right) \right\|_{C^m(M \times M)} \leq C \hbar^{N - 2n - 2m - 2k + 1}
\end{equation*}
for all $\hbar < \hbar_0$. Here, $\chi$ is a smooth cutoff function that is compactly supported in $M \bowtie M$ such that $\chi \equiv 1$ on a neighborhood of the diagonal in $M \times M$.
\end{theorem}

\begin{remark}
In fact, Bär and Pfäffle {\normalfont\cite{baerpfaeffle}} show Thm.\ \ref{ThmBP0} only for $W=0$. In this case, $\Phi_j=0$ whenever $j$ is odd. While this is not true for non-zero $W$, the transport equations \eqref{RecursiveTransportEquations} can be solved just as well and the proof of Thm.\ \ref{ThmBP} carries over basically without changes. 
\end{remark}


\section{Getzler Symbols} \label{SectionGetzlerSymbols}

This section is dedicated to the proof of the following theorem.

\begin{theorem} \label{CliffordSymbolOfHeatKernel}
Let $H_\hbar$ be the Witten operator of Example \ref{ExampleWittenOperator} acting on sections of $\Lambda \TMd$ and let $\mathrm{e}_\hbar \sum\nolimits_{j=0}^\infty \hbar^j \Phi_j$ be the corresponding formal heat kernel.
Then for each $j$, the coefficient $\Phi_j$ has Clifford bi-order at most $(j, j)$ and
\begin{equation}
  \sum_{j=0}^n \hbar^j \slashed{\sigma}_{j, j}(\Phi_j) 
  = \exp\bigl(  - t \,\slashed{\sigma}_{0, 0}(V) -  t \hbar \,\slashed{\sigma}_{1, 1}(W) - t\hbar^2 \slashed{\sigma}_{2, 2}(\mathbf{F}) \bigr).
\end{equation}
Here,
\begin{equation} \label{FatEff}
  \mathbf{F} =  -\frac{1}{8} R_{ijkl}\, \cc^i\cc^j \bb^k\bb^l.
\end{equation}
\end{theorem}

To establish a proof, we follow along the lines of Getzler's proof \citep{getzlerlocal} of the local Atiyah-Singer index theorem in the form explained in \citep{johnroe}: We define a symbol calculus on the space $\mathfrak{D}(M, \Lambda \TMd)$ of differential operators acting on sections of $\Lambda \TMd$ that takes into account the Clifford order. Then we show how to use it in the $\hbar$-dependent situation discussed here.

Following Roe, a symbol map is a homomorphism-like mapping from a filtered algebra to a graded algebra:

\begin{definition} \citep[Def.\ 12.5]{johnroe} \label{DefinitionSymbolMap}
Let $A$ be a filtered algebra and let $G$ be a graded algebra. A {\em symbol map} $\sigma_\bullet: A \longrightarrow G$ is a family of linear maps $A_k \longrightarrow G^k$ such that
\begin{enumerate}
  \item[(i)] If $a \in A_{l}$ with $l<k$, then $\sigma_k(a) = 0$
  \item[(ii)] If $a \in A_k$ and $b \in A_l$, then $\sigma_k(a)\sigma_l(b) = \sigma_{k+l}(ab)$.
\end{enumerate}
\end{definition}

We already introduced two symbol maps in Section \ref{SectionCliffordSymbol}. The Clifford symbol
\begin{equation*}
  \slashed{\sigma}_\bullet: \End(\Lambda V) \longrightarrow \Lambda V \,\stimes\, \Cl(-V)
\end{equation*}
and the Clifford bi-symbol. Here, $A = \End(\Lambda V)$ for a finite-dimensional vector space $V$ considered with either its $\Z$-filtration or its $\Z^2$-filtration. The corresponding graded algebras are $\Lambda V \,\stimes\, \Cl(- V)$ and $\Lambda V \stimes \Lambda V$, respectively (where the first is equipped with the grading that only takes into account the differential form degree).

Now we introduce Getzler's symbol. Let $(M, g)$ be a Riemannian manifold of dimension $n$. The algebra $\mathfrak{D}(M, \Lambda \TMd)$ of differential operators acting on sections of $\Lambda \TMd$ has a natural filtration by order, and comes with the map that assigns each operator its principal symbol. This is indeed a symbol map in the sense of Def.\ \ref{DefinitionSymbolMap}: The corresponding graded algebra is the algebra of sections of the bundle\footnote{Strictly speaking, this is an infinite-dimensional bundle. However, the homogeneous parts of this bundle are finite-dimensional, and we only consider homogeneous sections or finite linear combinations of these.}
\begin{equation*}
   \bigoplus_{k=0}^\infty S^k \TMd \otimes \Lambda \TMd \cong \Lambda \TMd[TM].
\end{equation*}
Here, $S^k \TMd$ denotes the $k$-th symmetric power of $\TMd$; this space is canonically isomorphic to the space of homogeneous polynomials of degree $k$ on $TM$ and the direct sum over these spaces for all $k \geq 0$ can be identified with a space of polynomials. In our case, we obtain the space of polynomials with coefficients in $\Lambda \TMd$, which we denote by $\Lambda \TMd[TM]$.

Now let us use all this data to define a new filtration on $\mathfrak{D}(M, \Lambda \TMd)$ that takes into account the Clifford order.

For $q \in M$, choose a chart $x$ around $q$ with $x(q)=0$. If we additionally trivialize the bundle $\Lambda \TMd$ by identifying fibers $\Lambda^k T^*_p\!M$ with $\Lambda^k T^*_q\!M$ for $p$ near $q$, any differential operator $P \in \mathfrak{D}(M, \Lambda \TMd)$ has a Taylor series with respect to these choices,
\begin{equation} \label{TaylorseriesOfP}
  P \sim \sum_{\alpha\beta} p_{\alpha\beta}\, x^\alpha \,\frac{\partial^{|\beta|}}{\partial x^\beta}, ~~~~~ \text{with}~~ p_{\alpha\beta} \in \End(\Lambda T^*_q\!M).
\end{equation}
Here, $\alpha$ and $\beta$ are multi-indices and we employ the usual conventions for these (as explained e.g.\ in \citep[p.\ 1]{shubin}).
Of course, this Taylor series depends heavily on the two choices made. Its order and the principal term, however, do not, allowing us to define the following.

\begin{definition}[$q$-Symbols] \label{DefinitionQSymbols}
Let $q \in M$. We say that $P$ is of $q$-order $k$ or less, if for each  $\alpha$ and $\beta$, $p_{\alpha\beta}$ has order less or equal to $k + |\alpha|-|\beta|$ in the Clifford filtration of $\End(\Lambda T^*_q\!M)$. In this case, its $k$-th {\em $q$-symbol} is
\begin{equation}
  \sigma_k^q(P) = \sum_{j =  k+|\alpha|- |\beta|-|\alpha|} \slashed{\sigma}_j(p_{\alpha\beta}) X^\alpha \frac{\partial^{|\beta|}}{\partial X^\beta}.
\end{equation}
Here, $X_j$ are the Euclidean coordinate functions on $T_qM$ induced by the chart $x$. $\sigma_k^q(P)$ is a differential operator on $T_qM$ with coefficients in the algebra $\mathcal{A}_q[T_qM]$, the space of $\mathcal{A}_q$-valued polynomials on $T_qM$, where
\begin{equation*}
\mathcal{A}_q := \Lambda T^*_q\!M \stimes \Cl(-T_qM).
\end{equation*}
We denote the space of such operators by $\mathfrak{P}(T_qM, \mathcal{A}_q)$. If an operator $P \in \mathfrak{D}(M, \Lambda \TMd)$ is of $q$-order $k$ or less for every $q \in M$, we can take its $q$-symbol at every point. This associates to $P$ its (global) {\em Getzler symbol} which we denote by $\sigma_k(P)$. This is a section of the bundle $\mathfrak{P}(TM, \mathcal{A})^k$.
\end{definition}

Under Euclidean coordinates $X_1, \dots, X_n$ on $T_qM$, we understand coordinates such that for all $v = v^j e_j \in T_qM$, we have $X^j(v) = v^j$ (here $e_1, \dots, e_n$ is an arbitrary but fixed orthonormal basis of $T_qM$). This implies $\frac{\partial}{\partial X_j}|_v = e_j$ for every $v \in T_qM$ under the canonical isomorphism $T_v T_qM \cong T_qM$. Each $X_j$ can be interpreted as a polynomial of degree one on $T_q M$.

The space $\mathfrak{P}(T_qM, \mathcal{A}_q)$ introduced in Def.\ \eqref{DefinitionQSymbols} is a graded algebra: We assigns degree one to differential operators $\frac{\partial}{\partial X^j}$ and one-forms $\eta \in T^*_q\!M \subseteq \Lambda T^*_q\!M$, and degree minus one to polynomials $X_j \in \R[T_qM]$. Elements of $\Cl(-T_qM)$ are assigned degree zero.

\begin{remark} \label{RemarkCliffordAndPSymbol}
The $q$-symbol is a refinement (and extension at the same time) of the Clifford filtration: If $P$ is an endomorphism (i.e.\ a differential operator of order zero in the usual filtration) and it has $q$-order $k$ or less (i.e.\ $\sigma^q_k(P) \in \mathcal{A}_q[T_qM]$), then the Clifford order of $P|_q$ is $k$ or less as well, and $\slashed{\sigma}_k(P|_q)$ is given by evaluating the polynomial $\sigma_k^q(P)$ at $X=0 \in T_qM$.
\end{remark}

\begin{proposition} \label{WelldefinednessOfSymbols}
The $q$-symbol is well-defined, i.e.\ its definition above is independent of the choices made. Furthermore, we have
\begin{equation} \label{MultiplicationRule}
  \sigma_{j+k}^q(P\circ Q) = \sigma_j^q(P)\circ\sigma_k^q(Q)
\end{equation}
whenever $P$ is of $q$-order $\leq j$ and $Q$ is of order $\leq k$.
\end{proposition}

\begin{proof}
If we have the well-definedness of the $q$-symbol map, \eqref{MultiplicationRule} follows directly from the the composition rules of differential operators in $\mathfrak{P}(T_qM, \mathcal{A}_q)$. Let us temporarily denote by $\sigma_k^{q, x}$ the $k$-th $q$-symbol with respect to a chart $x$. If we change from a chart $x$ to a chart $y$, then the charts have Taylor expansions with respect to each other, namely
\begin{align*}
    x &\sim Ay + \dots ~~~~~ \text{and} \\
    y &\sim Bx + \dots
\end{align*}
for some $A \in \mathrm{GL}(n)$ and $B = A^{-1}$, where the dots indicate terms of lower order in the $q$-filtration.
 Therefore
 \begin{equation*}
   \sigma_{|\alpha|}^{q, x}(x^\alpha) = X^\alpha = (AY)^\alpha = \sigma_{|\alpha|}^{q, y}((Ay)^\alpha) = \sigma_{|\alpha|}^{q, y}(x^\alpha),
 \end{equation*}
 because the order of $(Ay)^\alpha - x^{\alpha}$ is less than $|\alpha|$. A similar computation can be done for the differential operator $\partial^{|\beta|}/\partial x^\beta$, with the matrix $A$ replaced by $B$. 
 
Now under a trivialization, we understood a smooth identification of near-by fibers $\Lambda^k_p \TMd$ with $\Lambda^k T^*_q\!M$. A transformation $T$ from one such identification to another is of course the identity at $q$, so the $q$-order of $T - \mathrm{id}$ is negative which means that  the $q$-symbol is independent of the choice of trivialization. 
\end{proof}

\begin{lemma} 
Let $\nabla$ be the Levi-Civita connection on $\Lambda \TMd$ and choose Euclidean coordinates $X^1, \dots, X^n$ on $T_qM$. Then the operator $\nabla_i$ is of $q$-order $1$ and
\begin{equation} \label{SymbolOfCliffordConnection}
  \sigma_1^q(\nabla_i) = \frac{\partial}{\partial X^i} - \frac{1}{4}\sum\nolimits_{j=1}^n R_{ij} X^j
\end{equation}
where 
\begin{equation} \label{DefinitionRij}
  R_{ij} = \sum\nolimits_{k < l} R_{ijkl}\, \dd X^k \wedge \dd X^l
\end{equation}
are the $2$-form entries of the curvature tensor of $TM$ with respect to this basis.
\end{lemma}

\begin{proof}
Let $x$ be Riemannian normal coordinates about the point $q \in M$. From \cite[Lemma 4.14]{bgv}, we obtain that 
\begin{equation} \label{SymbolOfCliffordConnectionb}
  \nabla_i = \frac{\partial}{\partial x^i} - \frac{1}{8} R_{ijkl} x^j \cc^k \cc^l + \dots
\end{equation}
where the dots indicate terms that are of order zero or less in the $q$-filtration. Hence
  \begin{equation*}
   \sigma_1^q(\nabla_i) = \frac{\partial}{\partial X^i} - \frac{1}{8} R_{ijkl} X^j\, \dd X^k \wedge \dd X^l,
\end{equation*}
which is the claimed formula.
\end{proof}

\begin{remark}
  We adopt the convention $R_{ijkl} = \langle R(\partial_i, \partial_j)\partial_k, \partial_l\rangle$ which different from the one used in \citep{bgv}, also leading to a different sign of $R_{ij}$ (compare Lemma 4.14). Roe \cite{johnroe} uses the same sign convention for $R_{ijkl}$ as we do (compare Example 12.14) but chooses to define his $R_{ij}$ forms to be the negative of those defined in \eqref{DefinitionRij}. Our convention makes formula \eqref{SymbolOfCliffordConnectionb} coincide with \cite{getzlerlocal}. 
\end{remark}

\begin{remark}
By \citep[Thm.\ 12.13]{johnroe}, there is a unique symbol map 
\begin{equation*}
\sigma_\bullet: \mathfrak{D}(M, \Lambda \TMd) \longrightarrow \Gamma^\infty(M, \mathfrak{P}(TM, \mathcal{A})
\end{equation*} that satisfies formula \eqref{SymbolOfCliffordConnection} and restricts to the Clifford symbol defined in Section \ref{SectionCliffordSymbol} when restricted to an operator of order zero (i.e.\ a section of $\End(\Lambda\TMd)$). Therefore, our global Getzler symbol coincides with the one defined in \citep[Thm.\ 12.13]{johnroe}, restricted to the special case of the Hodge Laplacian considered here (the general definition works on any Clifford bundle). In Getzler's original paper \cite{getzlerlocal}, a similar calculus is implicit, but the notion of a symbol map is not formalized.
\end{remark}

Getzler's main observation was that the $q$-symbol of a Dirac operator is a harmonic oscillator on $T_qM$, meaning the following:

\begin{proposition} {\normalfont\citep[Prop.\ 12.17]{johnroe}} \label{PropSymbolOfD2}
Let $D = \dd + \dd^*$. For every $q \in M$, $D^2$ has $q$-order $2$ and its $q$-symbol with respect to orthogonal coordinates on $T_q M$ is
\begin{equation} \label{SymbolOfD2}
  \sigma_2^q(D^2) = - \sum\nolimits_{i=1}^n \left( \frac{\partial}{\partial X^i} - \frac{1}{4}\sum\nolimits_{j=1}^n R_{ij} X^j \right)^2 + \slashed{\sigma}_{2}(\mathbf{F}) \in \mathfrak{P}(T_qM, \mathcal{A}_q),
\end{equation}
where $\mathbf{F}$ was defined in \eqref{FatEff}.
\end{proposition}

\begin{proof}
This follows at once from the Weizenböck formula for the Euler operator $D = \dd + \dd^*$, which is
\begin{equation} \label{Weizenboeck}
  \dd \dd^* + \dd^* \dd = D^2 = \nabla^* \nabla + \frac{1}{4}\mathrm{scal} + \mathbf{F},
\end{equation}
(compare (3.16) and p.\ 149 in \cite{bgv}) and formula \eqref{SymbolOfCliffordConnection}).
\end{proof}

The $q$-symbols of the other terms in the Witten operator $H_\hbar = \hbar^2 D^2 + \hbar W + V$ are straightforward to calculate; they are given by
\begin{align}
 \label{SymbolOfW} {\sigma}_1^q(W) &= \slashed{\sigma}_{1, \bullet}(W) = \sum_{i=1}^n \dd x^i \stimes \bb(\nabla_i \xi)\\
 \label{SymbolOfV}  {\sigma}_0^q(V) &= \slashed{\sigma}_{0, \bullet}(V) = |\xi|^2.
\end{align}

We now need to extend the symbol calculus to the coefficients $\Phi_j$ of the formal heat kernel of $H_\hbar$. We abbreviate $\Eh := \Lambda \TMd$. For each  $q \in M$ fixed and each $t>0$, $\Phi_j(\,\cdot\,, q, t)$ is a section of the vector bundle $\Eh^* \otimes \Eh_q$ over $M$. Therefore, with respect to a chart $x$ around $q$ and a trivialization that identifies fibers of $\Eh_p$ with $\Eh_q$ for $p$ near $q$, it has a Taylor series
\begin{equation*}
  \Phi_j(\,\cdot\,, q, t) \sim \sum_\alpha \Phi_{j\alpha}(q, t) \, x^\alpha, ~~~ \text{where} ~~~ \Phi_{j\alpha}(q, t)  \in \Eh^*_q \otimes \Eh_q \cong \End(\Eh_q),
\end{equation*}
where the coefficients depend smoothly on $q$ and $t$. We say that $\Phi_j$ has $q$-order $k$ or less if $\Phi_{j\alpha}(q, t) \in \End(\Eh_q)$ has Clifford order less or equal to $k + |\alpha|$ for each $t$, and in that case, we define its $k$-th $q$-symbol as
\begin{equation*}
  \sigma_k^q(\Phi_j) = \sum_{m + |\alpha| = k} \slashed{\sigma}_m\bigl(\Phi_{j\alpha}(q, t)\bigr) X^\alpha \in \mathcal{A}[T_pM].
\end{equation*}
Note that this also depends on $t$ in addition to the usual dependence on $q$.

The well-definedness can be shown as in the proof of Prop.\ \ref{WelldefinednessOfSymbols} and a formal  computation with Taylor series shows that if the $q$-order of $\Phi_j$ is $l$ or less, we have the multiplication property
\begin{equation} \label{MultiplicationProperty}
  \sigma_k^q(P) \sigma_l^q(\Phi_j) = \sigma_{k+l}^q(P \,\Phi_j)
\end{equation}
for any differential operator $P \in \mathfrak{D}(M, \Lambda \TMd)$ of order at most $k$, which is supposed to act on the first entry of $\Phi_j$. 

\begin{theorem} \label{FormulaForHeatSymbol}
 For each $j = 0, 1, 2, \dots$, and each $q \in M$,   $\Phi_j$ is of $q$-order at most $j$ and for the "heat symbol" $\sigma^q(k_\hbar) := \mathbf{e}_\hbar\sum_{j=0}^n \hbar^j \sigma_j^q(\Phi_j)$, we have the formula 
  \begin{equation} \label{FormulaForSymbolOfK}
    \sigma^q(k_\hbar) = u(X, R, t\hbar^2) \exp\bigl(- t\, \slashed{\sigma}_0(V) - t \hbar\,\slashed{\sigma}_1(W) -t\hbar^2 \slashed{\sigma}_2(\mathbf{F}) \bigr)
  \end{equation}
where
\begin{equation*}
  u(X, R, t) = (4 \pi t)^{-n/2} \hat{A}(t R) 
  \exp \left( -\frac{1}{4t}\left\langle X, \frac{tR}{2}\coth \left(\frac{tR}{2}\right)X \right\rangle \right)
\end{equation*}
is the Mehler kernel with the $\hat{A}$-form
\begin{equation*}
 \hat{A}(R) = \det\nolimits^{1/2} \left( \frac{R/2}{\sinh (R/2)} \right).
\end{equation*}
\end{theorem}

Let us first see why this implies Thm.\ \ref{CliffordSymbolOfHeatKernel}.

\begin{proof}[of Thm.\ \ref{CliffordSymbolOfHeatKernel}]
With a view on Remark \ref{RemarkCliffordAndPSymbol}, we have
\begin{equation*}
  \sum_{j=0}^n \hbar^j \slashed{\sigma}_{j, \bullet}(\Phi_j) \!\!= (4\pi t\hbar^2)^{n/2}\sigma^q(k_\hbar)|_{X=0} = \hat{A}(t \hbar^2 \!R)\exp\bigl(- t\, \slashed{\sigma}_{0, \bullet}(V) - t \hbar\,\slashed{\sigma}_{1, \bullet}(W) -t\hbar^2 \slashed{\sigma}_{2, \bullet}(\mathbf{F}) \bigr).
\end{equation*}
Furthermore $\hat{A}(t\hbar^2R) = 1 + \sum_{j=1}^{n} \hbar^{j} \hat{A}_{j}(t)$ for some $j$-form\footnote{In fact, $\hat{A}_{j}(t)=0$ unless $j$ is divisible by four, as $\hat{A}$ is an even function of $R$.} $\hat{A}_{j}(t) \in \Omega^{j}(M)$ and the $\exp$-part is of the form $\sum_{j=0}^n \hbar^j \slashed{\sigma}_{j, \bullet}(E_j(t))$ for some $E_j(t) \in \End(\Lambda T^*_q M)$ of order $(j, j)$ in the Clifford bi-filtration because the exponents have bi-order $(0, 0)$, $(1, 1)$ and $(2, 2)$, respectively and come with appropriate powers of $\hbar$. Therefore
\begin{equation*}
  \sum_{j=0}^n \hbar^j \slashed{\sigma}_{j, \bullet}(\Phi_j) = \sum_{j=0}^n \hbar^j \sum_{k+l=j} \hat{A}_k(t) \slashed{\sigma}_{l, \bullet}(E_l(t))
\end{equation*}
Because $\hat{A}_k(t)$ has order $(k, 0)$ in the Clifford bi-filtration, we now find that $\slashed{\sigma}_{j, j}(\Phi_j) = \hat{A}_0(t) \slashed{\sigma}_{j, j}(E_j(t)) = \slashed{\sigma}_{j, j}(E_j(t))$. This gives the theorem.
\end{proof}

We now proceed with the proof of Thm.\ \ref{FormulaForHeatSymbol}. For each $q \in M$, the "symbolic operator"
\begin{equation}
  \sigma^q(H_\hbar) := \hbar^2 \sigma_2^q(D^2) + \hbar \sigma_1^q(W) + \sigma_0^q(V)
\end{equation}
is of Schrödinger type as an operator on $T_qM$, acting on $C^\infty$-functions with values in $\mathcal{A}_q$ (i.e.\ sections of a trivial vector bundle). Therefore, we can consider the "symbolic heat equation"
\begin{equation} \label{HeatEquationOnTqM}
  \bigl( \partial_t + \sigma^q(H_\hbar) \bigr) \mathbf{k}_\hbar = 0
\end{equation}
on $T_qM$.
We use bold letters for all objects associated to this symbolic equation to distinguish them from those associated to the equation down on $M$. By Thm.\ \ref{ThmBP0}, there exists a unique formal heat kernel
\begin{equation*}
  \hat{\mathbf{k}}_\hbar(X, Y, t) = \mathbf{e}_\hbar(X-Y, t)\sum_{j=0}^\infty \hbar^j\, \boldsymbol{\Phi}_j(X, Y, t),
\end{equation*}
where the $\boldsymbol{\Phi}_j$ solve the transport equations corresponding to the heat equation \eqref{HeatEquationOnTqM} and 
\begin{equation*}
  \mathbf{e}_\hbar(X, t) = (4\pi t\hbar^2)^{-n/2} \exp\left( -\frac{1}{4t\hbar^2} |X|^2 \right).
\end{equation*}
is the Euclidean fundamental solution on $T_q M$.

Looking at formula \eqref{SymbolOfD2}, the unique connection $\boldsymbol{\nabla}$ on $T_q M$ such that $\sigma_2^q(D^2) - \boldsymbol{\nabla}^*\boldsymbol{\nabla}$ is of order zero (compare \cite[Prop.\ 2.5]{bgv})  is clearly given by
\begin{equation*}
  \boldsymbol{\nabla}_i = \frac{\partial}{\partial X^i} - \sum\nolimits_{j=1}^n R_{ij} X^j = \sigma_1^q(\nabla_i).
\end{equation*}
Because $T_qM$ is flat as a Riemannian manifold, $\mathbf{G}\equiv 0$ and $\boldsymbol{\mathcal{V}} = \sum_j X^j \frac{\partial}{\partial X^j}$, hence
\begin{equation*}
  \boldsymbol{\nabla}_{\boldsymbol{\mathcal{V}}} = \sum\nolimits_{j=1}^n X^j \frac{\partial}{\partial X^j} - \sum\nolimits_{ij=1}^n R_{ij} X^i X^j = \sigma^q_0(\nabla_{\mathcal{V}})
\end{equation*}
The recursive transport equations \eqref{RecursiveTransportEquations} for the operator $\sigma^q(H)$ therefore can be written as
\begin{equation} \label{TransportEquationOnTqM}
  \left( t\frac{\partial}{\partial t} + \sigma^q_0(\nabla_{\mathcal{V}}) + t \sigma^q_0(V) \right) \boldsymbol{\Phi}_j(\,\cdot\,, Y) = - t \sigma^q_1(W)\boldsymbol{\Phi}_{j-1}(\,\cdot\,, Y) - t\sigma^q_2(D^2)\boldsymbol{\Phi}_{j-2}(\,\cdot\,, Y).
\end{equation}

\begin{lemma}
  For each $j = 0, 1, 2, \dots$ and each $q \in M$, the coefficient $\Phi_j$ is of $q$-order at most $j$. 
\end{lemma}

\begin{proof}
 By definition, for each $q \in M$, the $\Phi_j$ satisfy the transport equations
  \begin{equation} \label{TransportEquation2}
    \left( t \frac{\partial}{\partial t} + \nabla_{\mathcal{V}} + G + tV\right) \Phi_j(\,\cdot\,, q) = - t W \Phi_{j-1}(\,\cdot\,, q) - tL \Phi_{j-2}(\,\cdot, q)
  \end{equation}
  with initial condition $\Phi_0(q, q, 0) = \mathrm{id}_{\Eh_q}$. Straightforward calculation shows that $G(\,\cdot\,, q)$ vanishes at $q$, hence $\sigma^q_0(G) = 0$. Let $\Phi_0$ be of $q$-order $k$. Taking the $k$-th symbol on both sides shows (by multiplicativity \eqref{MultiplicationProperty}) that $\sigma^q_k(\Phi_0)$ solves
\begin{equation*}
  \left( t \frac{\partial}{\partial t} + \sigma^q_0(\nabla_{\mathcal{V}}) + t \sigma^q_0(V) \right)\sigma^q_k(\Phi_0) = 0,
\end{equation*}
which is just the first transport equation on $T_qM$ \eqref{TransportEquationOnTqM}. By Thm.\ \ref{ThmBP0}, there is a unique solution for each initial value. Because $\Phi_0(q, q, 0) = \mathrm{id}_{\Eh_p}$, we have 
\begin{equation*}
  \sigma^q_k(\Phi_0)|_{(X, t)=0} = \slashed{\sigma}_k(\Phi_0(q, q, 0)) = \slashed{\sigma}_k(\mathrm{id}_{\Eh_q}) = 0 ~~~~~ \text{if}~~ k >0,
\end{equation*}
so the initial value is zero and $\sigma^q_k(\Phi_0)$ vanishes for all $X$ and $t$. This shows that $\Phi_0$ has $q$-order zero for all $t$.
Now by induction and multiplicativity of symbols, the right hand side of the $j$-th transport equation \eqref{TransportEquation2} has $q$-order $\leq j$. Suppose that $\Phi_j$ has $q$ order $k>j$. Taking the $k$-th $q$-symbol on both sides shows that $\sigma^q_k(\Phi_j)$ solves
\begin{equation*}
 \left( t \frac{\partial}{\partial t} + \sigma^q_0(\nabla_\mathcal{V}) + t \sigma^q_0(V)\right)\sigma^q_k(\Phi_j) = 0.
\end{equation*}
Again, this is a transport equation. By Thm.\ \ref{ThmBP0}, there is a unique solution for each right hand side, and zero is a solution, hence $\sigma^q_k(\Phi_j) = 0$ whenever $k>j$. But this means that $\Phi_j$ is of $q$-order $\leq j$.
\end{proof}

This gives the following  corollary.

\begin{corollary} \label{CorollaryHeatSymbol}
  The terms $\sigma^q_j(\Phi_j)$ solve the recursive transport equations \eqref{TransportEquationOnTqM} for $Y = 0$ and we have
  \begin{equation} \label{Equation1}
   \boldsymbol{\Phi}_j|_{Y=0} = \sigma^q_j(\Phi_j).
   \end{equation}
Furthermore, because $\sigma^q_j(\Phi_j) = 0$ whenever $j>n$, the formal heat kernel $\hat{\mathbf{k}}_\hbar$ is actually a finite sum and we have
\begin{equation*}
    \sigma^q(k_\hbar) = \hat{\mathbf{k}}_\hbar|_{Y=0}.
\end{equation*}
\end{corollary}

\begin{proof}
Taking the $j$-th $q$-symbol on both sides of the $j$-th equation \eqref{TransportEquation2} gives exactly the transport equations \eqref{TransportEquationOnTqM} for $\sigma^q(H_\hbar)$. Equation \eqref{Equation1} follows from the uniqueness statement of Thm.\ \ref{ThmBP0}.
\end{proof}

Let us now finish the proof of the Thm.\ \ref{FormulaForHeatSymbol}.

\begin{proof}[of Thm.\ \ref{FormulaForHeatSymbol}]
First, one verifies that the right hand side of \eqref{FormulaForSymbolOfK} is a solution to the heat equation \eqref{HeatEquationOnTqM}. The result that $u(X, R, t)$ satisfies
\begin{equation*}
  \left(\frac{\partial}{\partial t}  - \sum\nolimits_{i=1}^n \left( \frac{\partial}{\partial X^i} - \frac{1}{4}\sum\nolimits_{j=1}^n R_{ij} X^j \right)^2\right)u(X, R, t) = 0, ~~~~~~~~ u_0 = \delta_0 \cdot \mathrm{id}_{\Eh_q}
\end{equation*}
is usually called Mehler's formula, see \citep[chapter 4.2]{bgv}. Substitution $t \mapsto \hbar^2t$ and straightforward calculation then shows that the right hand side of \eqref{FormulaForSymbolOfK} solves \eqref{HeatEquationOnTqM}.

Now explicitly expanding the Taylor series, one verifies that
\begin{equation*}
  u(X, R, t\hbar^2) \exp\bigl(- t\, \slashed{\sigma}_0(V) - t \hbar\,\slashed{\sigma}_1(W) -t\hbar^2 \slashed{\sigma}_2(\mathbf{F}) \bigr)
  = \mathbf{e}_\hbar(X, t) \boldsymbol{\Phi}_\hbar(X, t),
\end{equation*}
for some power series $\boldsymbol{\Phi}_\hbar = \sum_{j=0}^\infty \hbar^j \boldsymbol{\Phi}_j$ with coefficients in $\mathcal{A}[T_q M]$, each $\boldsymbol{\Phi}_j$ being a polynomial in both $X$ and $t$. As seen in section \ref{SectionHeatkernel}, the $\boldsymbol{\Phi}_j$ have to fulfill the recursive transport equations \eqref{TransportEquationOnTqM}, and by uniqueness (Thm.\ \ref{ThmBP0}), we get
\begin{equation*}
  \boldsymbol{\Phi}_\hbar = \sum_{j=0}^n \hbar^j \sigma^q_j(\Phi_j)
\end{equation*}
as by Corollary \ref{CorollaryHeatSymbol}, the $\sigma^q_j(\Phi_j)$ as well solve these transport equations.
\end{proof}


\section{The McKean-Singer Formula and its Consequences} \label{SectionMcKeanSingr}

From now on, let $(M, g)$ be an $n$-dimensional closed \mbox{Riemannian} manifold. Again, let $D_\hbar$ be Witten's perturbed Euler operator of Example \ref{ExampleWittenOperator}. The McKean-Singer formula \citep[Thm. 3.50]{bgv} states that for all $t>0$, we have
\begin{equation*}
  \ind(D_\hbar) = \int_M \str k_\hbar(q, q, t)\, \dd q.
\end{equation*}
where we integrate over the volume density associated to the Riemannian metric of $M$.
The index of the Euler operator $D = \dd + \dd^*$ on the exterior algebra equipped with the grading given by $\Xi$ is well-known to be equal to the Euler characteristic $\chi(M)$. On the other hand, the index is a topological invariant, i.e.\ it is the same for every Dirac operator on a given vector bundle \citep[Thm.\ 3.51]{bgv}. Therefore, the index of $D_\hbar$ is also equal to $\chi(M)$. Expanding $k_\hbar$ in its semiclassical expansion and using that the asyptotics are uniform over $M$ by Thm.\ \ref{ThmBP}, we get
\begin{equation*}
  \chi(M) \sim_{\hbar \searrow 0} (4\pi t \hbar^2)^{-n/2} \sum_{j=0}^\infty \hbar^j \int_M \str \Phi_j(q, q, t)\, \dd q
\end{equation*}
for any $t>0$ fixed. The left hand side of this equation is independent of $\hbar$, so by uniqueness of asymptotic expansions, all coefficients on the right-hand-side except the term constant in $\hbar$ must vanish. Hence "$\sim$" must in fact be an equality and we get
\begin{equation}
    \chi (M) = (4 \pi t)^{-n/2} \int_M \str \Phi_n(q, q, t)\, \dd q ~~~~ \text{for all} ~~~~ t>0.
\end{equation}
By Corollary \ref{SupertraceFormula}, the supertrace of $\Phi_n(t)$ can be calculated in terms of the Clifford bi-symbol via the formula
\begin{align*} 
    \str\Phi_n (q, q, t) &=(-1)^{\frac{n(n+1)}{2}}\, 2^n \,\bigl\langle \slashed{\sigma}_{n,n}\bigl(\Phi_n(q, q, t)\bigr), \vol \, \widehat{\otimes} \, \vol \bigr\rangle \\
    &= (-1)^{\frac{n(n+1)}{2}}\, 2^n \sum_{j=0}^n \,\bigl\langle \slashed{\sigma}_{j,j}\bigl(\Phi_n(q, q, t)\bigr), \vol \, \widehat{\otimes} \, \vol \bigr\rangle
\end{align*}
where the last equality holds because the pairing of $\slashed{\sigma}_{j, j}(\Phi_n)$ with $\vol \stimes \vol$ vanishes if $j < n$.
On the other hand, by Thm.\ \ref{CliffordSymbolOfHeatKernel} we have
\begin{equation*}
\sum_{j=0}^n \slashed{\sigma}_{j,j}\bigl(\Phi_j(q, q, t)\bigr) = \exp\bigl(- t |\xi|^2 - t \,\slashed{\sigma}_{1, 1}(W) -t\slashed{\sigma}_{2,2}(\mathbf{F}) \bigr).
\end{equation*}
where we used $\slashed{\sigma}_{0,0}(V) = |\xi|^2$. 
Therefore, we get that for all $t>0$, we have
  \begin{equation} \label{IntegralFormula1}
    \chi (M) = (-1)^{\frac{n(n+1)}{2}}(\pi t)^{-n/2} \int_M \Bigl\langle \exp\bigl(- t |\xi|^2 - t \,\slashed{\sigma}_{1,1}(W) -t \slashed{\sigma}_{2,2}(\mathbf{F}) \bigr), \vol \stimes \vol \Bigr\rangle.
  \end{equation}

\begin{remark}
  The integrand above is in fact the pullback of Mathai and Quillen's  Thom form \cite{mathaiquillen}
\begin{equation*} 
  U = (2\pi)^{-n/2} T\bigl( \exp(- |\boldsymbol{\eta}|^2/2 + i \nabla \boldsymbol{\eta} + F ) \bigr) \in \Omega^n(\TMd)
\end{equation*}  
along the section $\xi_t := (2t)^{1/2} \xi$. 
Regarding the terms appearing above, $\boldsymbol{\eta}$ is the tautological section in $\Gamma^\infty(\TMd, \pi^* \TMd)$ that maps $\eta \mapsto \eta$, $T$ is the Berezin integral on the second component, $F$ is the Riemann tensor considered as $(2, 2)$-form on $\TMd$ and $\nabla$ is the Levi-Civita connection on $\TMd$, both pulled back to $\TMd$ via the canonical projection $\pi$ (see \citep[Chapter 1.6]{bgv} or \citep[Chapter 3]{zhang}).

Note that $\slashed{\sigma}_{1,1}(W)$ is not quite equal to $\nabla \boldsymbol{\eta}$ because of the appearance of a super tensor product, whence the lack of the factor $i$ in \eqref{IntegralFormula1}. By one possible definition, the Euler form is the pullback $\iota^*U$ along the inclusion of $\iota:M \longrightarrow T^*M$ as the zero section, which corresponds to $t=0$. Therefore, \eqref{IntegralFormula1} is exactly the interpolation formula
\begin{equation*}
  \chi(M) = \int_M \xi_t^* U,
\end{equation*}
which can be found in the proof of \citep[Thm.\ 1.56]{bgv}.

In this sense, formula \eqref{IntegralFormula1} can be seen as a generalization of a result of Mathai \cite{mathai}, who showed how to get the Thom form on the tangent bundle from the heat kernel of the Laplacian.
\end{remark}

Let us now evaluate \eqref{IntegralFormula1} without knowing anything about Thom forms. Note that $\slashed{\sigma}_{2, 2}(\mathbf{F})$, $\slashed{\sigma}_{1, 1}(W)$ and $|\xi|^2$ all commute so that
\begin{equation*}
  \exp\bigl( - t \slashed{\sigma}_{2, 2}(\mathbf{F}) - t  \slashed{\sigma}_{1, 1}(W) - t |\xi|^2 \bigr) = e^{-t|\xi|^2} \sum_{k=0}^n (-t)^k \sum_{j=0}^k \frac{\slashed{\sigma}_{2, 2}(\mathbf{F})^j \,\slashed{\sigma}_{1, 1}(W)^{k-j}}{j!(k-j)!}
\end{equation*}
We are interested in the $(n, n)$-form part, which is produced by the summands above where $n = 2j + (k-j) = k+j$. Therefore, we define the functions $\alpha_k$ by
\begin{equation} \label{FormulaAj}
  \alpha_k := \frac{(-1)^k}{(n-k)!(2k-n)!} \Bigl\langle \slashed{\sigma}_{2, 2}(\mathbf{F})^{n-k} \slashed{\sigma}_{1, 1}(W)^{2k-n}, \vol \stimes \vol \Bigr\rangle
\end{equation}
whenever $2k\geq n$. Then
\begin{equation*}
  \Bigl\langle \exp\bigl( - t \slashed{\sigma}_{2, 2}(\mathbf{F}) - t  \slashed{\sigma}_{1, 1}(W) - t |\xi|^2 \bigr), \vol \stimes \vol \Bigr\rangle = e^{-t|\xi|^2}\sum_{2k \geq n}^{n} t^k \alpha_k.
\end{equation*}
Using formula \eqref{IntegralFormula1}, we get the following formula for the index density.

\begin{theorem} \label{MainTheoremWithAj}
  Let $M$ be an $n$-dimensional closed \mbox{Riemannian} manifold. Then for all $t>0$, we have
  \begin{equation} \label{IntegralFormulaAj}
    \chi (M) = (-1)^{\frac{n(n+1)}{2}}(\pi t)^{-n/2} \sum_{2k \geq n}^{n} t^k \int_M \alpha_k \, e^{-t|\xi|^2}
  \end{equation}
  where the functions $\alpha_k$ are defined in \eqref{FormulaAj}.
\end{theorem}

\begin{proposition}
If $n$ is even, then the function $\alpha_{n/2}$ is given by
\begin{equation*} 
\alpha_{n/2} = \frac{1}{8^{n/2} (n/2)!} \sum_{\tau, \sigma \in \mathcal{S}_n} \mathrm{sgn}(\tau)\mathrm{sgn}(\sigma) R_{\tau(1)\tau(2)\sigma(1)\sigma(2)} \cdots R_{\tau(n-1)\tau(n)\sigma(n-1)\sigma(n)}.
\end{equation*}
\end{proposition}

\begin{proof}
  By formula \eqref{FormulaAj},
\begin{equation*}
  \alpha_{n/2} = \frac{(-1)^{n/2}}{(n/2)!} \, \Bigl\langle\slashed{\sigma}_{2, 2}(\mathbf{F})^{n/2}, \vol \stimes \vol\Bigr\rangle,
\end{equation*}
so the proposition follows directly from the formula
\begin{equation} \label{CliffordSymbolF}
\sigma_{2, 2}(\mathbf{F}) = -\frac{1}{8}\sum_{i, j, k, l = 1}^n R_{ijkl}\, e^ie^j\stimes e^k e^l,
\end{equation}
compare \eqref{FatEff}.
\end{proof}

The term
\begin{equation} \label{LipschitzKilling}
  \Omega = \frac{(-1)^{n/2}}{2^n (n/2)!} \sum_{\tau, \sigma \in \mathcal{S}_n} \mathrm{sgn}(\tau)\mathrm{sgn}(\sigma) R_{\tau(1)\tau(2)\sigma(1)\sigma(2)} \cdots R_{\tau(n-1)\tau(n)\sigma(n-1)\sigma(n)}
\end{equation}
 is called Killing-Lipschitz curvature (or $n$-th order sectional curvature) and up to a factor coincides with the Pfaffian of the curvature tensor (compare \citep[(3.39)]{zhang}). Taking the limit $t \downarrow 0$ in \eqref{IntegralFormulaAj}, we obtain the following classical theorem (note that if $n$ is odd, there is no term of order zero in $t$ in \eqref{IntegralFormulaAj} so that the integrand vanishes in the limit).

\begin{theorem}[Gauss-Bonnet-Chern] {\normalfont \citep[Thm.\ 1]{cherncurvature}}\label{GaussBonnetChern}
  Let $M$ be an $n$-dimensional closed \mbox{Riemannian} manifold. If $n$ is even, then its Euler characteristic is given by the integral formula
  \begin{equation}
    \chi(M) = (2\pi)^{-n/2}\int_M \Omega,
  \end{equation}
  where $\Omega$ is the Lipschitz-Killing curvature as in \eqref{LipschitzKilling}. If $n$ is odd, then $\chi(M) = 0$.
\end{theorem}

\begin{proposition}
We have
\begin{equation} \label{FormulaForAlphan}
  \alpha_{n} = (-1)^{\frac{n(n+1)}{2}} \det\nolimits_{g} \bigr(\nabla \xi\bigr)
\end{equation}
where the determinant of the $(0, 2)$-tensor $\nabla \xi$ is calculated with help of the metric.
\end{proposition}

\begin{proof}
By \eqref{FormulaAj}, we have 
\begin{equation*}
\alpha_{n}  = \frac{(-1)^n}{n!} \Bigl\langle \slashed{\sigma}_{1, 1}(W)^n, \vol \stimes \vol \Bigr\rangle.
\end{equation*}
We have
\begin{equation*}
  \slashed{\sigma}_{1, 1}(W) = \sum_{i=1}^n \dd x^i \stimes \nabla_i \xi =\sum_{ij} (\nabla_i \xi)_j \dd x^i \stimes \dd x^j
\end{equation*}
By multiplication rules of the super tensor product, we generally have
  \begin{equation} \label{MultiplicationRules}
    \left( \sum_{\alpha, \beta =a}^b c_{\alpha\beta} \,\dd x^\alpha \stimes \dd x^\beta\right)^k
    \!\!=(-1)^{\frac{k(k-1)}{2}}\!\!\!\!\!\!\!\sum_{\substack{\alpha_1, \dots, \alpha_k=a \\ \beta_1, \dots, \beta_k =a}}^b  \!\!\!\!\!c_{\alpha_1\beta_1} \cdots  c_{\alpha_k\beta_k}\, \dd x^{\alpha_1} \wedge \dots \wedge \dd x^{\alpha_k} \stimes \dd x^{\beta_1} \wedge \dots \wedge \dd x^{\beta_k}.
  \end{equation}
Hence
\begin{align*}
  \Bigl\langle \slashed{\sigma}_{1, 1}(W)^n, \vol \stimes \vol \Bigr\rangle  
  &= (-1)^{\frac{n(n-1)}{2}} \sum_{\substack{\alpha_1, \dots, \alpha_k=1 \\ \beta_1, \dots, \beta_k =1}}^n  \!\!\!\!\!(\nabla_{\alpha_1} \xi)_{\beta_1} \cdots  (\nabla_{\alpha_n} \xi)_{\beta_n} \\
   &= (-1)^{\frac{n(n-1)}{2}} \sum_{\tau, \sigma \in \mathcal{S}_n} \mathrm{sgn} (\sigma)\mathrm{sgn}(\tau)\, (\nabla_{\tau(1)} \xi)_{\sigma(1)} \cdots (\nabla_{\tau(n)} \xi)_{\sigma(n)} \\
    &= (-1)^{\frac{n(n-1)}{2}} n! \det(\nabla \xi),
\end{align*}
by the Leibnitz formula for determinants. The result follows
\end{proof}

Now we can evaluate the integral \eqref{IntegralFormulaAj} with the method of stationary phase, where we take the limit $t \uparrow \infty$ in \eqref{IntegralFormulaAj}. This allows us to express the Euler characteristic as a sum of integrals over the set of  critical points of $\xi$. The task is to explicitly calculate certain jets of the functions $\alpha_j$, which appear in the stationary phase expansion.

We always assume that the zero locus of $\xi$ is a disjoint union of finitely many submanifolds of $M$. The easiest situation is the case that all critical submanifolds are zero-dimensional, which is dealt with in the following corollary. The degenerate case is dealt with in the next section.

\begin{theorem}[Poincaré-Hopf] \label{ThmPoincareHopf}
Let $M$ be a closed Riemannian manifold. Suppose that the vector field $X$ has only non-degenerate critical points, i.e.\ whenever we have $X(p) = 0$, then $\det(\nabla X|_p) \neq 0$. Then
\begin{equation*}
  \chi(M) = \sum_{\{X(p) = 0\}} (-1)^{\nu(p)},
\end{equation*}
where $\nu(p)$ is the number of negative eigenvalues of $\nabla X|_p$.
\end{theorem}

This is a special case of Thm.\ \ref{ThmMorseBott}. We give a separate proof, as it is quite short and hopefully indicates how the proof works in the general case. Note that we do not require $M$ to be oriented (compare e.g.\ \citep[Thm.\ 1.58]{bgv}). 

\begin{proof}
Set $\xi := X^\flat$, the adjoint form to $X$ given by the metric (i.e.\ "lowering the indices"). Notice that the set $C$ of points where $|X|^2 = |\xi|^2 =0$ is just the set $\mathrm{Crit}_\phi$ of critical points of $\phi := |X|^2$. The method of stationary phase (see e.g.\ \citep[Prop.\ 5.2]{dimassisjoestrand} or compare Lemma \ref{LemmaStationaryPhase} below) states that if $\phi$ fulfills $\phi > 0$ on the compact support of a function $\alpha$ except at the point $p$, where $\phi(p)=0$ and $\det \nabla^2\phi > 0$, then we have
\begin{equation*}
  \lim_{t \uparrow \infty} t^{k} \int_M \alpha \, e^{-t \phi} 
  = \begin{cases} 
  (2\pi)^{n/2} \det \bigl(\nabla^2 \phi|_p\bigr)^{-1/2} \alpha(p) &\text{if}~~k = n/2 \\
  0 & \text{if}~~ k < n/2
  \end{cases}.
\end{equation*}
Here, $\nabla^2 \phi|_p$ is the Hessian of $\phi$ at $p$. By a partition of unity argument, taking the limit $t \uparrow \infty$ in \eqref{IntegralFormulaAj} yields
\begin{align*}
  \chi(M) = \lim_{t \uparrow \infty}\, \left(\frac{t}{\pi}\right)^{n/2} \int_M \alpha_{n} \,e^{-t|X|^2}
  = 2^{n/2}\!\!\sum_{\{ X(p) = 0\}} \!\!  \det \bigl( \nabla^2 |X|^2|_p\bigr)^{-1/2} \alpha_n(p)
\end{align*}
At $p$, we have $\nabla^2_{Y, Z} |X|^2 = \langle \nabla_Y X, \nabla_Z X\rangle$
so that the determinant is given by $\det\nolimits( \nabla^2 |X|^2) = 2^n \det(\nabla X)^2$. Therefore, using \eqref{FormulaForAlphan}
\begin{equation} \label{DeterminantOfX2}
\det \bigl(\nabla^2 |X|^2\bigr)^{-1/2} \alpha_n(p) = \frac{\det(\nabla X)}{2^{n/2} \bigl|\det(\nabla X)\bigr|} = \frac{(-1)^{\nu(p)}}{2^{n/2}}
\end{equation}
and the theorem follows.
\end{proof}


\section{The Degenerate Case} \label{SectionDegenerateCase}

In this section, we generalize Thm.\ \ref{ThmPoincareHopf} to the case that we have critical submanifolds instead of critical points.

\begin{definition} \cite{bottnondegenerate} \label{DefMorseBott}
  Let $M$ be a manifold. A function $\phi$ is called {\em Morse-Bott}, if
  \begin{enumerate}
  \item[(i)] The set of critical points $\mathrm{Crit}_\phi = \{ p \in M \mid \dd \phi|_p = 0 \}$ is a disjoint union of finitely many submanifolds of $M$.
  \item[(ii)] For each connected component $C \subseteq \mathrm{Crit}_\phi$, the Hessian $\nabla^2 \phi$ of $X$ is non-degenerate when restricted to the normal bundle $N C$.
\end{enumerate}
The second condition {\em a priori} depends on the choice of a Riemannian metric, but turns out to be independent of this choice. 
\end{definition}

\begin{definition}
  The {\em index} $\nu(C)$ of a connected critical submanifold $C \subseteq \mathrm{Crit}_\phi$ of $M$ is the dimension of the biggest subspace $V$ of $T_pM$ such that the bilinear form $\nabla^2\phi|_p$ restricted to $V$ is negative definite, where $p$ is some point in $C$.
\end{definition}

This is independent of the point $p \in C$ because $\nu(C)$ is locally constant, as follows from the Morse-Bott Lemma (see for example \cite{morsebott}).

\begin{theorem}[Degenerate Poincaré-Hopf] \label{ThmMorseBott}
Let $M$ be an closed manifold. Suppose that $\phi$ is Morse-Bott as defined above and let $\mathrm{Crit}_\phi = C_1 \coprod \cdots \coprod C_k$ for connected submanifolds $C_j$. Then 
  \begin{equation*}
    \chi(M) = \sum_{j=1}^k (-1)^{\nu(C_j)} \chi(C_j).
  \end{equation*}
\end{theorem}

Usually one proves this theorem using the degenerate Morse inequalities \citep[Thm.\ 2.14]{bismut}. We will instead start from Thm.\ \ref{MainTheoremWithAj} and evaluate the integrand in the limit $t \uparrow \infty$ with the method of stationary phase. To this end, we first make some general observations regarding integrals of the form
\begin{equation*}
  I(t, \alpha) = \int_M e^{-t |\dd \phi|^2} \alpha.
\end{equation*}
Here $M$ is an $n$-dimensional Riemannian manifold, $\phi$ is a Morse-Bott function, and $\alpha \in C^\infty_c(M)$. Clearly, $I(t, \alpha)$ is a smooth function on $(0, \infty)$.

\begin{lemma} \label{LemmaStationaryPhase}
Assume that $\mathrm{Crit}_\phi$ consists of exactly one submanifold $C$, which has co-dimension $m$. Then 
\begin{equation*}
  \lim_{t \uparrow \infty} t^{m/2} I(t, \alpha) = \left(\frac{\pi}{2}\right)^{m/2}\,\int_C \frac{\alpha}{\,\bigl|\det(\nabla^2\phi|_{NC})\bigr|},
\end{equation*}
where we integrate over the canonical measure on $C$ induced by the Riemannian metric.
\end{lemma}

\begin{proof}
If $\phi$ is Morse-Bott, then $|\dd \phi|^2$ is a Morse-Bott function as well.
Assume that there exists a diffeomorphism $\kappa= (\kappa_1, \dots, \kappa_m, \widetilde{\kappa}): M \supseteq U \longrightarrow (V\subseteq \R^m) \times C$ (where $U$ is a tubular neighborhood of $C$ and $V$ is a neighborhood of zero in $\R^m$) such that
\begin{equation*}
  |\dd \phi|^2\bigr|_U = \kappa_1^2 + \dots + \kappa_m^2
\end{equation*}
Then, assuming that $\alpha$ is supported in $U$, we have
\begin{equation}
  t^{-m/2} I(t, \alpha) = \int_C \left(\int_V e^{-t{|v|^2}} (\alpha\circ \kappa^{-1})(x, p) \,G(x, p) \,\dd v\right) \dd p,
\end{equation}
where $G = |\det \dd \kappa|^{-1}$ is the Jacobian determinant coming from the transformation formula. It is well-known that $t^m e^{-t|v|^2}$ converges as a distribution to $\pi^{m/2} \delta_0$ in the limit $t \uparrow \infty$, where $\delta_0$ is the delta distribution at zero in $\R^m$. Straightforward calculation (e.g.\ in local coordinates) yields
\begin{equation*}
  G(0, p) = \det \bigl(\nabla^2 |\dd \phi|^2\bigr|_{N_pC}\bigr)^{-1/2} = 2^{-m/2} \bigl|\det\bigl(\nabla^2\phi|_{N_pM}\bigr)\bigr|^{-1}.
\end{equation*}
This proves the lemma for $\alpha$ compactly supported in $U$, and under the assumption of the existence of a diffeomorphism $\kappa$ as above.

For general $\alpha$, we can use a partition of unity to split it up into a function $\alpha_1$ that is compactly supported in $U$ and a function $\alpha_2$ with compact support disjoint from $C$. For the first, we can use the argument above, and the second part is of order $t^{-\infty}$, because $|\dd \phi|^2 \geq \varepsilon$ on the support of $\alpha_2$, for some $\varepsilon > 0$.

Regarding the diffeomorphism $\kappa$, such a diffeomorphism exists at least locally on $C$ by \cite{morsebott}, and the statement again follows by splitting up $\alpha$ with a partition of unity.
\end{proof}

\begin{proposition} \label{PropStationaryPhase}
Under the assumptions of Lemma \ref{LemmaStationaryPhase}, $I(t, \alpha)$ has a complete asymptotic expansion as $t \uparrow \infty$, namely
\begin{equation} \label{StationaryPhaseExpansion}
  I(t, \alpha) \sim \left(\frac{\pi}{2t}\right)^{m/2}\sum_{j=0}^\infty t^{-j} \int_C \frac{L^j \alpha}{j!\,\bigl|\det(\nabla^2 \phi|_{NC})\bigr|}.
\end{equation}
Here, $L$ is a second order differential operator defined on a neighborhood of $C$ that for each $q \in C$ has the $q$-symbol
\begin{equation*}
  \sigma_2^q(L) = \sum_{\alpha, \beta=1}^mA_q^{\alpha\beta}\frac{\partial^2}{\partial X^\alpha \partial X^\beta}~~~~~~   \text{with}~~~~~~A_q =
  \frac{1}{4}
  \begin{pmatrix}
    \bigl( \nabla^2\phi|_{N_qC}\bigr)^{-2}& 0 \\
    0 & 0
  \end{pmatrix}
\end{equation*}
in the decomposition $T_qM = N_qC \oplus T_qC $.
\end{proposition}

The asymptotic expansion \eqref{StationaryPhaseExpansion} is often called stationary phase expansion, though one usually considers {\em imaginary} emponents.

\begin{proof}
It is more convenient to consider $J(s, \alpha):= s^{-m/2}I(s^{-1}, \alpha)$ instead. The idea is then to construct an operator $L$ defined on some neighborhood of $C$ such that
\begin{equation} \label{HeatEqCritical}
  \frac{\partial}{\partial s} J(s, \alpha) = J(s, L\alpha).
\end{equation}
If we make the ansatz
\begin{equation*}
  L^*\alpha:= \div \bigl( A \cdot \grad \alpha \bigr) + \bigl\langle v, \grad \alpha \bigr\rangle,
\end{equation*}
with an endomorphism field $A$ and a vector field $v$, calculate both sides of \eqref{HeatEqCritical} and order by powers of $s$, we obtain that the endomorphism $A$ and the vector field $v$ have to fulfill the {\em coefficient equations}
\begin{align}
 \label{CoefficientEquation1} |\dd \phi|^2 &= \bigl\langle \grad|\dd \phi|^2, A \cdot \grad |\dd \phi|^2 \bigr\rangle \\
 \label{CoefficientEquation2} \bigl\langle v, \grad |\dd \phi|^2 \bigr\rangle &=  \frac{m}{2} - \div \bigl( A \cdot \grad |\dd \phi|^2 \bigr)
\end{align}
We have
\begin{equation*}
  \grad |\dd \phi|^2 = 2 \langle \nabla \dd \phi, \dd \phi \rangle^\sharp = 2 \nabla^2 \phi \cdot \grad \phi.
\end{equation*}
so that \eqref{CoefficientEquation1} becomes
\begin{equation*}
  |\grad \phi|^2 = 4 \bigl\langle \grad \phi, \nabla^2\phi \cdot A \cdot \nabla^2\phi \cdot \grad \phi \bigr\rangle
\end{equation*}
Therefore, if we find a vector field $V$ such that
\begin{equation}  \label{TransportEquationV}
  \nabla_{\grad \phi} V = 2\nabla^2\phi \cdot \grad \phi,
\end{equation}
then $A = (\nabla V)^*\nabla V$ solves \eqref{CoefficientEquation1}. Equation \eqref{TransportEquationV} is a transport equation that is singular at $C$. Because the right-hand side vanishes on $C$ as well, it admits a unique smooth solution for all given smooth initial values on $C$ (this is explained in \citep[Thm.\ 2.3, Thm.\ 7.1]{ludewigsource}), so in particular there exists a unique endomorphism field $A$ of the form $A = (\nabla V)^*\nabla V$, defined on some neighborhood of $C$, such that the induces operator $L$ has the claimed principal symbol.

Having constructed such an $A$, we find $\div\bigl(A \cdot \grad |\dd\phi|^2\bigr) \equiv \frac{m}{2}$ on $C$. This means that the right-hand-side of the second coefficient equation \eqref{CoefficientEquation2} vanishes on $C$, and again, the equation has a solution of the form $v = \grad f$, defined on some neighborhood of $C$. 

Together, this shows that on some neighborhood of $C$, there exists an operator $L$ with the claimed principal symbol such that \eqref{HeatEqCritical} holds. Therefore, we can expand in a Taylor series around any $\varepsilon >0$.
\begin{equation*}
  J(s, \alpha) = \sum_{j=0}^N \frac{(s-\varepsilon)^j}{j!} J(\varepsilon, L^j\alpha) + \frac{1}{N!}\int_{\varepsilon}^s (s-r)^N J(r, L^{N+1}\alpha) \dd r
\end{equation*}
By \eqref{LemmaStationaryPhase}, we can take the limit $\varepsilon\downarrow 0$ in this formula, which implies that we have the claimed asymptotic expansion, at least if the support of $\alpha$ is compactly contained in a neighborhood of $C$ where $L$ is defined. The general case follows with a partition of unity argument, as explained before.
\end{proof}

We now apply the stationary phase expansion \eqref{StationaryPhaseExpansion} to the integral formula \eqref{IntegralFormulaAj} for the Euler characteristic, where we set $\xi := \dd \phi$. We may and will henceforth assume that $\mathrm{Crit}_\phi$ consists of exactly one connected submanifold $C$ of co-dimension $m$; the general case follows once more from a partition of unity argument. This gives
\begin{align*}
  \chi(M) &= (-1)^{\frac{n(n+1)}{2}}(\pi t)^{-n/2} \sum_{2k \geq n}^n t^{k} I(t, \alpha_k) \\
  &\sim (-1)^{\frac{n(n+1)}{2}}(\pi t)^{-n/2} \left(\frac{\pi}{2t}\right)^{m/2} \sum_{k \geq n}^n \sum_{j=0}^\infty t^{k - j} \!\!\! \int_C \frac{L^j\alpha_k}{j! \bigl|\det(\nabla^2 \phi|_{NC})\bigr|}.
\end{align*}
Now we invoke the same argument as earlier: Because the left-hand-side is independent of $t$, all coefficients of the asymptotic expansion of the right-hand-side must be zero, except the constant term, and we must have equality in fact. The zero-order terms are those where $k=n/2 + m/2 + j$ ($k$ and $j$ being integers) so that we have zero-order terms if and only if $m$ and $n$ have the same parity, or equivalently, if $C$ is even-dimensional. We have obtained the following partial result.

\begin{lemma}
If the critical set of $\phi$ consists of exactly one submanifold $C$ of co-dimension $m$, we have the following integral formula
\begin{equation} \label{StationaryPhaseFormula}
  \chi(M) = (-1)^{\frac{n(n+1)}{2}}\pi^{m/2-n/2}\sum_{j=0}^{n/2-m/2} \int_C \frac{L^j\alpha_{n/2 + m/2 + j}}{j! \bigl|\det(\nabla^2\phi|_{NC})\bigr|}.
\end{equation} 
if and only if $C$ is even-dimensional. Otherwise, $\chi(M)=0$.
\end{lemma}

\begin{remark}
This implies already that for any Morse-Bott function $\phi$ on $M$, $\mathrm{Crit}_\phi$ can be a union of non-empty odd-dimensional manifolds only if $\chi(M)=0$.
\end{remark}

We can therefore assume that $C$ is even-dimensional, i.e.\ that $m$ and $n$ have the same parity. The task is now to explicitly calculate the integrands $L^j \alpha_{n/2+m/2+j}$ on $C$. It will turn out that $\alpha_{n/2+m/2+j}$ has $q$-order $-2j$ for every $q \in C$. Therefore,
\begin{equation} \label{SymbolEvaluation}
  L^j\alpha_{n/2+m/2+j}(q) = \sigma_2^q(L)^j \sigma_{-2j}^q(\alpha_{n/2+m/2+j})\bigr|_{X=0}
\end{equation}
where $\sigma^q_\bullet$ contains the $q$-symbol of Section \ref{SectionGetzlerSymbols}.

To make our calculations easier, we fix a point $p \in C$ and work in a special chart around the point $p$, constructed as follows: Choose an orthonormal basis $\nu_1|_p, \dots, \nu_m|_p$ of $N_pC$ and obtain an orthonormal frame $\nu_1, \dots, \nu_m$ of $NC$ over a neighborhood of $p$ by parallel translation along geodesics in $C$ emanating from $p$. Furthermore, choose Riemannian normal coordinates $\widetilde{x}$ in $C$ centered around the point $C$. Now define the chart $x: M \supset U \longrightarrow V \subset \R^m \times \R^{n-m}$ by
\begin{equation*}
  x^{-1}(v, w) := \exp_{q}\Bigl(\sum\nolimits_{\alpha=1}^m v^\alpha \nu_\alpha|_q\Bigr), ~~~~~ q = \widetilde{x}^{-1}(w),
\end{equation*}
where $U$ is a suitable neighborhood of $p$ in $M$. Write $\partial_\alpha$ for the directional derivatives with respect to this chart. Then we have the following standard result.

\begin{lemma}
For $q \in C\cap U$, we have $\partial_\alpha|_q = \nu_\alpha|_q$ if $\alpha \leq m$ and $\partial_\alpha|_q \in T_qC$ if $\alpha > m$. The Christoffel symbols at the point $p$ are
\begin{equation}
\Gamma_{\alpha\beta}^\gamma|_p = 
\begin{cases} \mathrm{II}_{\alpha\beta}^\gamma|_p ~~~~~ &\text{if}~~ \gamma \leq m ~\text{and}~ \alpha, \beta > m \\
0 & \text{otherwise} \\
\end{cases}.
\end{equation}
Here $\mathrm{II}_{\alpha\beta}^\gamma$ are the coefficients of the second fundamental form $\mathrm{II}|_p: T_p C \times T_pC \longrightarrow N_p C$, compare {\normalfont\citep[II.2.1]{chavel}}.
Furthermore, in the decomposition $T_pM =N_pC \oplus T_pC$,
\begin{equation} \label{HessianOfPhi}
  \nabla^2 \phi|_p ~\hat{=}~
  \begin{pmatrix}
    \phi_{\alpha\beta} & 0\\ 0 & 0
  \end{pmatrix}.
\end{equation}
with respect to the basis $\partial_1|_p, \dots, \partial_n|_p$, where $\phi_{\alpha\beta} = \partial_\alpha\partial_\beta \phi|_p$.
\end{lemma}

\begin{lemma}
For every $p \in C$ and $0 \leq j \leq n/2 - m/2$, the function $\alpha_{n/2+m/2+j}$ has $p$-order $-2j$ or less, and its $p$-symbol is given by
\begin{equation*}
  \sigma_{-2j}^p(\alpha_{n/2 + m/2 +j}) = {c}_j \,\det \bigl(\nabla^2 \phi|_{N_pC}\bigr)  
 \!\!\! \sum_{\substack{\sigma, \tau \in \mathcal{S}_{m, n}}} \!\!\!
    \mathrm{sgn}(\tau)\,\mathrm{sgn}(\sigma)~ \Theta_{\tau\sigma}(X)
\end{equation*}
where ${c}_j$ is a dimensional constant and the $\Theta_{\tau\sigma}(X)$ are homogeneous polynomials of degree $2j$ on $T_pM$ that only depend on the curvature of $M$, the second fundamental form of $C$ and the Hessian of $\phi$ at $p$; see  \eqref{DefinitionOfTheta} and \eqref{TheConstantCj}. Here, the summation is over all $\sigma, \tau$ in the group $\mathcal{S}_{m, n}$ of permutations of the numbers $\{m+1, \dots, n\}$.
\end{lemma}

\begin{proof}
  By \eqref{SymbolOfW}, we have 
  \begin{equation*}
  \slashed{\sigma}_{1,1}(W) =
  \sum_{\alpha, \beta=1}^n \nabla_{\alpha\beta}^2\phi \,\dd x^\alpha \stimes \dd x^\beta =
   \!\sum_{\alpha, \beta=1}^n \!\left(\partial_\alpha \partial_\beta \phi - \sum_{\gamma=1}^n\Gamma_{\alpha\beta}^\gamma \partial_\gamma \phi \right) \dd x^\alpha \stimes \dd x^\beta.
  \end{equation*}
   By \eqref{HessianOfPhi}, $\nabla^2_{\alpha\beta} \phi$ has $p$-order zero or less with $\sigma_0(\nabla^2_{\alpha\beta} \phi) = \phi_{\alpha\beta}$ if $\alpha, \beta \leq m$ and if $\alpha, \beta > m$, then $\nabla^2_{\alpha\beta} \phi$ has $p$-order $-1$ or less with
  \begin{equation} \label{SymbolOfNablaAlphaBeta}
    \sigma_{-1}^p\bigl(\nabla^2_{\alpha\beta} \phi\bigr) =  
    -\sum_{\gamma, \delta =1}^m \mathrm{II}_{\alpha\beta}^\gamma\, \phi_{\gamma\delta} X^\delta.
\end{equation}
  Therefore,
  \begin{align*}
    \slashed{\sigma}_{1,1}(W)^{m+2j} 
    &= \binom{m+2j}{m} \underbrace{ \left( \sum_{\alpha, \beta =1}^m \nabla^2_{\alpha\beta}\phi \, \dd x^\alpha \stimes \dd x^\beta\right)^{\!\!\!m}\!\!}_{:= \mathcal{W}_1}
   \underbrace{ \left( \sum_{\alpha, \beta =m+1}^n \!\!\!\!\nabla^2_{\alpha\beta}\phi \, \dd x^\alpha \stimes \dd x^\beta \right)^{\!\!\!2j}\!\!\!}_{:= \mathcal{W}_2} \!+ O(|x|^{2j+1})
  \end{align*}
  The pre-factor accommodates for the fact that we selected $m$ factors out of $m+2j$ factors to split $\slashed{\sigma}_{1,1}(W)^{m+2j}$ up this way. 
  For the first factor we then get (compare \eqref{MultiplicationRules})
  \begin{align*}
    \mathcal{W}_1|_p &= (-1)^{\frac{m(m-1)}{2}}\!\!\!\!\!\!\sum_{\substack{\alpha_1, \dots, \alpha_m=1 \\ \beta_1, \dots, \beta_m =1}}^m  \!\!\!\!\!\!\!\phi_{\alpha_1\beta_1} \cdots  \phi_{\alpha_m\beta_m}\, \dd x^{\alpha_1} \wedge \dots \wedge \dd x^{\alpha_m} \stimes \dd x^{\beta_1} \wedge \dots \wedge \dd x^{\beta_m} \\
    &= (-1)^{\frac{m(m-1)}{2}}\!\!\!\!\sum_{\tau, \sigma \in \mathcal{S}_m} \,\!\!\!\!\!\mathrm{sgn}(\tau)\mathrm{sgn}(\sigma)\, \phi_{\tau(1)\sigma(1)} \cdots  \phi_{\tau(m)\sigma(m)} \, \dd x^1 \!\wedge \dots \wedge \dd x^m \stimes \dd x^1 \!\wedge \dots \wedge \dd x^m \\
    &= (-1)^{\frac{m(m-1)}{2}} m! \,\det \bigl(\nabla^2 \phi|_{N_pC}\bigr) \,\dd x^1 \!\wedge \dots \wedge \dd x^m \stimes \dd x^1 \!\wedge \dots \wedge \dd x^m
  \end{align*}
  and for the second factor, we have
  \begin{align*}
    \mathcal{W}_2 &= \left( -\!\!\!\!\!\sum_{\alpha, \beta =m+1}^n \nabla^2_{\alpha\beta}\phi\,\, \dd x^\alpha \stimes \dd x^\beta \right)^{2j} 
   = (-1)^j \sum_{A, B} \Phi_{AB} \,\dd x^A \stimes \dd x^B,
  \end{align*}
  where $A = (\alpha_1, \dots, \alpha_{2j})$ and $B= (\beta_1, \dots, \beta_{2j})$ run over all $2j$-tuples of numbers between $m+1$ and $n$. By the considerations above, $\Phi_{AB}$ is of $p$-order $-2j$ and its $p$-symbol is
  \begin{equation*}
    \sigma_{-2j}^p\bigl(\Phi_{AB}\bigr) = \sum_{\substack{\gamma_1, \dots, \gamma_{2j} =1 \\ \delta_1, \dots, \delta_{2j}=1}}^m 
    \mathrm{II}_{\alpha_1\beta_1}^{\gamma_1} \cdots \mathrm{II}_{\alpha_{2j}\beta_{2j}}^{\gamma_{2j}} \phi_{\gamma_1\delta_1}  \cdots \phi_{\gamma_{2j}\delta_{2j}} \,\, X^{\delta_1} \cdots X^{\delta_{2j}},
  \end{equation*}
compare \eqref{SymbolOfNablaAlphaBeta}. Now remember that
\begin{equation*}
  \alpha_{n/2+m/2 +j} = \frac{(-1)^{n/2 + m/2+j}}{(n/2-m/2-j)!(m+2j)!} \, \Bigl \langle \slashed{\sigma}_{2, 2}(\mathbf{F})^{n/2-m/2-j} \slashed{\sigma}_{1, 1}(W)^{m+2j} \Bigr\rangle.
\end{equation*}
By the calculations above, we obtain
\begin{align*}
  \sigma_{-2j}^p(\alpha_{n/2 + m/2 +j}) &= {c}_j \,\det \bigl(\nabla^2 \phi|_{N_pC}\bigr)  
 \!\!\!\!\! \sum_{\substack{\sigma, \tau \in \mathcal{S}_{m, n}}} \!\!\!\!\!
    \mathrm{sgn}(\tau)\mathrm{sgn}(\sigma)~ \Theta_{\tau\sigma}(X)
\end{align*}  
for a dimensional constant ${c}_j$, where $\Theta_{\tau\sigma}(X)$ is the homogeneous polynomial
  \begin{equation} \label{DefinitionOfTheta}
   \Theta_{\tau\sigma}(X) = \sigma_{-2j}^p\bigl(\Phi_{AB}\bigr) 
    \!\!\!\!\!\prod_{s=m/2+j+1}^{n/2}\!\!\!\!\! R_{\tau(2s-1)\tau(2s)\sigma(2s-1)\sigma(2s)}
    ~~~~\text{for}~~~~ \substack{A=(\tau(m+1), \dots, \tau(m+2j))\\B = (\sigma(m+1), \dots, \sigma(m+2j))}.
  \end{equation}
  Here we used formula \eqref{CliffordSymbolF} to calculate $\slashed{\sigma}_{2,2}(\mathbf{F})^{n/2-m/2-j}$. In the case that $n$ (and hence also $m$) is odd, the above product is supposed to run over half-integer numbers $s$.
  
For the constant ${c}_j$, we find
\begin{equation}
\begin{split} \label{TheConstantCj}
  {c}_j &= \frac{(-1)^{n/2 + m/2+j}}{(n/2-m/2-j)!(m+2j)!} \binom{m+2j}{m} (-1)^{\frac{m(m-1)}{2} +j}m!  \left(-\frac{1}{8}\right)^{n/2-m/2-j}   \\
  &= \frac{(-1)^{n+\frac{m(m-1)}{2}}}{8^{n/2-m/2}}\frac{(-1)^j8^j}{(n/2-m/2-j)!(2j)!}
  \end{split}
\end{equation}
which only depends on $j$, $m$ and $n$.
\end{proof}

\begin{lemma}
For each $p \in C$, we have
\begin{equation}
(-1)^{\nu(C)} \,\Omega_C|_p = (-1)^{\frac{n(n+1)}{2}}{2^{n/2-m/2}} \sum_{j=0}^{n/2-m/2} \frac{L^j\alpha_{n/2 + m/2 + j}(p)}{j! \bigl|\det(\nabla^2\phi|_{N_pC})\bigr|}
\end{equation}
where $\Omega_C$ is the Lipschitz-Killing curvature of $C$, as in \eqref{LipschitzKilling}.
\end{lemma}

This proves Thm.\ \ref{ThmMorseBott} using the stationary phase formula \eqref{StationaryPhaseFormula} and the Gauss-Bonnet-Chern Theorem \ref{GaussBonnetChern}.

\begin{proof}

  Denote by $\phi^{\alpha\beta}$ the entries of the inverse matrix of $\bigl(\phi_{\alpha\beta}\bigr)_{\alpha\beta \leq m}$. Then the operator $\sigma_2^p(L)$ from Prop.\ \ref{PropStationaryPhase} is given by
  \begin{equation*}
    \sigma_2^p(L) = \frac{1}{4}\sum_{\alpha, \beta, \gamma=1}^m \phi^{\alpha\gamma}\phi^{\beta\gamma} \frac{\partial^2}{\partial X^\alpha \partial X^\beta}.
  \end{equation*}
  To calculate $\sigma_2^p(L)^j \sigma_{-2j}^p\bigl(\Phi_{AB}\bigr)$, we need to split the $2j$ indices into groups of two, for which there are $2^{-j}(2j)!$ possibilities. Therefore
  \begin{align*}
    \sigma_2^p(L)^j \sigma_{-2j}^p\bigl(\Phi_{AB}\bigr) 
   &=\frac{(2j)!}{2^j}\prod_{s=1}^j\sigma_2^p(L)\left(\sum_{\substack{\gamma \gamma^\prime \delta \delta^\prime =1}}^m \mathrm{II}_{\alpha_{2s-1}\beta_{2s-1}}^{\gamma}\mathrm{II}_{\alpha_{2s}\beta_{2s}}^{\gamma^\prime} \phi_{\gamma\delta}\phi_{\gamma^\prime\delta^\prime} X^\delta X^{\delta^\prime} \right) \\
    &= \frac{(2j)!}{2^j4^j}\prod_{s=1}^j\left(\sum_{\substack{\gamma =1}}^m \mathrm{II}_{\alpha_{2s-1}\beta_{2s-1}}^{\gamma}\mathrm{II}_{\alpha_{2s}\beta_{2s}}^{\gamma} + \sum_{\substack{\gamma^\prime =1}}^m \mathrm{II}_{\alpha_{2s-1}\beta_{2s-1}}^{\gamma^\prime}\mathrm{II}_{\alpha_{2s}\beta_{2s}}^{\gamma^\prime} \right)\\
    &= \frac{(2j)!}{8^{j}} 2^j\prod_{s=1}^j \bigl\langle \mathrm{II}_{\alpha_{2s-1}\beta_{2s-1}},  \mathrm{II}_{\alpha_{2s}\beta_{2s}} \bigr\rangle.
  \end{align*}
 Put together and using \eqref{SymbolEvaluation} now gives
  \begin{align*}
    L^j\alpha_{n/2 + m/2 +j} &= {c}_j \,\det \bigl(\nabla^2 \phi|_{N_pC}\bigr)  
 \!\!\! \sum_{\substack{\sigma, \tau \in \mathcal{S}_{m, n}}} \!\!\!
    \mathrm{sgn}(\tau)\,\mathrm{sgn}(\sigma)~ \sigma_2^p(L)^j\Theta_{\tau\sigma}(X)\bigr|_{X=0} \\
    &= \frac{(-1)^{n+\frac{m(m-1)}{2}}}{8^{n/2-m/2}}\frac{(-1)^j8^j}{(n/2-m/2-j)!(2j)!} \,\det \bigl(\nabla^2 \phi|_{N_pC}\bigr) \frac{(2j)!}{8^j} \Upsilon_j \\
    &= \frac{2^{n/2}(-1)^{n + \frac{m(m-1)}{2}}}{8^{n/2 - m/2}(n/2-m/2-j)!}\,\det \bigl(\nabla^2 \phi|_{N_pC}\bigr) (-1)^j \Upsilon_j
  \end{align*}
  where
  \begin{align*}
   \Upsilon_j &= 2^j\!\!\!\!\sum_{\sigma, \tau \in \mathcal{S}_{m,n}}  \!\!\!\!\mathrm{sgn}(\tau)\mathrm{sgn}(\sigma)\!\!
   \prod_{s=m/2}^{m/2+j} \bigl\langle \mathrm{II}_{\tau(2s-1)\sigma(2s-1)},  \mathrm{II}_{\tau(2s)\sigma(2s)} \bigr\rangle 
    \!\!\!\!\!\!\!\prod_{s=m/2+j+1}^{n/2}\!\!\!\!\!\!\! R_{\tau(2s-1)\tau(2s)\sigma(2s-1)\sigma(2s)} \\
    &= \!\!\sum_{\sigma, \tau \in \mathcal{S}_{m,n}}  \!\!
    \mathrm{sgn}(\tau)\mathrm{sgn}(\sigma)
   \prod_{s=m/2}^{m/2+j} S_{\tau(2s-1)\tau(2s)\sigma(2s-1)\sigma(2s)}
    \!\!\!\!\!\prod_{s=m/2+j+1}^{n/2}\!\!\!\!\! R_{\tau(2s-1)\tau(2s)\sigma(2s-1)\sigma(2s)}
  \end{align*}
  with 
  \begin{equation*}
    S_{\alpha\beta\gamma\delta} = \bigl\langle \mathrm{II}_{\alpha\gamma}, \mathrm{II}_{\beta\delta}\bigr\rangle - \bigl\langle \mathrm{II}_{\beta\gamma}, \mathrm{II}_{\alpha\delta}\bigr\rangle
  \end{equation*}

The Gauss formula \citep[Thm.\ II.2.1]{chavel} states that the entries $\widetilde{R}_{\alpha\beta\gamma\delta}$ of the curvature tensor of $C$ are given by the formula
\begin{equation*}
  \widetilde{R}_{\alpha\beta\gamma\delta} = R_{\alpha\beta\gamma\delta} - S_{\alpha\beta\gamma\delta}.
\end{equation*}
Therefore,
\begin{equation} \label{UpsilonSum}
 \sum_{\sigma, \tau \in \mathcal{S}_{m,n}}  \!\!\!\!\mathrm{sgn}(\tau)\mathrm{sgn}(\sigma) \!\!\!\prod_{s=1}^{n/2-m/2}\!\!\! \widetilde{R}_{\tau(2s-1)\tau(2s)\sigma(2s-1)\sigma(2s)} = \!\!\!\sum_{j=0}^{n/2-m/2}\!\!\!(-1)^j\binom{n/2-m/2}{j} \Upsilon_j.
\end{equation}
In total, we obtain
\begin{align*}
{2^{-m/2}} \!\!\!\sum_{j=0}^{n/2-m/2} \!\!\!\frac{L^j\alpha_{n/2 + m/2 + j}(p)}{j! \bigl|\det(\nabla^2\phi|_{N_pC})\bigr|} 
&= \!\!\!\sum_{j=0}^{n/2-m/2}\!\!\! \frac{\det \bigl(\nabla^2 \phi|_{N_pC}\bigr)}{j! \bigl|\det(\nabla^2\phi|_{N_pC})\bigr|} \frac{2^{n/2-m/2}(-1)^{n + \frac{m(m-1)}{2}}}{8^{n/2 - m/2}(n/2-m/2-j)!}\, (-1)^j\Upsilon_j \\
&= (-1)^{\nu(C)}\frac{(-1)^{n + \frac{m(m-1)}{2}}}{2^{n-m}(n/2-m/2)!} \sum_{j=0}^{n/2-m/2} (-1)^j \binom{n/2-m/2}{j} \Upsilon_j
\end{align*}
Finally, we use \eqref{UpsilonSum} and the fact that
\begin{equation*}
  n+ \frac{n(n+1)}{2} + \frac{m(m-1)}{2} = \frac{n}{2} - \frac{m}{2} \mod 2
\end{equation*}
because $n$ and $m$ have the same parity to obtain the stated result.
\end{proof}

\begin{remark}
It is a curious observation that the coefficients $\alpha_k$ vanish exactly to the order such that only there top-order terms matter in the stationary phase evaluation. In particular, this also means that only the principal symbol of the stationary phase operator $L$ contributes to the end result.
\end{remark}

\bibliography{literature}

\end{document}